\documentclass[11pt]{article}
\usepackage[latin1]{inputenc}
\usepackage[english]{babel}
\usepackage{amssymb}
\usepackage{amsmath,amsthm}
\usepackage[normalem]{ulem}

\newtheorem{defn}{Definition}[section]
\newtheorem{thm}{Theorem}[section]
\newtheorem{prop}{Proposition}[section]
\newtheorem{cor}{Corollary}[section]
\newtheorem{lem}{Lemma}[section]
\newtheorem{rmk}{Remark}[section]
\newtheorem{exmp}{Example}[section]

\makeatletter

\let \al=\alpha
\let \be=\beta
\let \var=\varphi
\let \vare=\varepsilon

\let \de=\delta
\let \th=\theta
\let \la=\lambda

\let \ga=\gamma
\let \p=\partial
\let \q=\quad
\let \qq=\qquad
\let \med=\medskip
\let \smal=\smallskip
\let \dps=\displaystyle
\let \ul=\underline

\let \ul=\underline
\let \ol=\overline

\newcommand{\R}{\mathbb{R}}
\newcommand{\N}{\mathbb{N}}

\begin{document}

\vspace{0.5cm}
{\bf {\Large \centerline{ Asymptotic behaviour for a class of non-monotone }

\centerline{delay differential systems with applications}}}

 \

 \centerline{\scshape Teresa Faria\footnote{Corresponding author.}
\footnote{
 Departamento de Matem\'atica and CMAF-CIO, Faculdade de Ci\^encias,\\
 Universidade de Lisboa, Campo Grande, 1749-016 Lisboa, Portugal\\
teresa.faria@fc.ul.pt}, Rafael Obaya\footnote{
Departamento de Matem\'{a}tica
Aplicada, E. Ingenier\'{\i}as Industriales,
 and  IMUVA,\\
 Instituto de
Matem\'{a}ticas, Universidad de Valladolid,
47011 Valladolid, Spain\\
rafoba@wmatem.eis.uva.es}, Ana M. Sanz
\footnote{
Departamento de Did\'{a}ctica de las Ciencias Experimentales, Sociales y de la Matem\'{a}tica, \\
  Facultad de Educaci\'{o}n, and  IMUVA,\\ Instituto de
Matem\'{a}ticas, Universidad de Valladolid,
 34004 Palencia, Spain\\
anasan@wmatem.eis.uva.es}}

\

\

{\it Suggested running head:} Asymptotic behaviour for a class of non-monotone DDEs

\

\

\begin{abstract}
The paper concerns a class of $n$-dimensional non-autonomous delay differential equations  obtained by adding a non-monotone delayed perturbation to a linear homogeneous cooperative system of ordinary differential equations. This family  covers a wide set of models used  in structured population dynamics. By exploiting the stability and the monotone character of the linear ODE, we establish sufficient conditions for both the extinction of all  the populations and the permanence of the system. In the case of DDEs with autonomous coefficients (but possible time-varying delays), sharp results are obtained, even in the case of a reducible community matrix. As a sub-product,
 our results  improve some criteria  for autonomous systems published in recent literature. As an important illustration,
 the extinction, persistence and permanence of a non-autonomous Nicholson system with patch structure and multiple time-dependent delays are analysed.
\end{abstract}


\med



\section{Introduction}
\setcounter{equation}{0}

 In the last decades,  delay differential equations  (DDEs) with patch structure have been largely employed in population dynamics and other fields, since by  capturing several features of a  heterogeneous environment, they may  provide quite realistic models. Structured systems of differential equations have been used in population models when the populations  are distributed over   different classes (e.g. due to age, size or   different  food-rich patches), in disease models with several compartments for a host population, in leukemia models where the  infected cells may become inactive, and in a variety of other situations where the  transition among the several classes should be considered. See e.g. \cite{Cushing,  MetzDiek,Smith11,SmithThi}. Naturally, time delays should be incorporated in such systems to express the maturation time of biological species, the incubation period of diseases, the maturation time of blood cells and several  other attributes.


The paper is concerned with a  family of non-autonomous  DDEs
written in abstract form as
\begin{equation}\label{i1}
x'(t)=A(t)x(t)+f(t,x_t),\q t\ge 0
\end{equation}
where $A(t)$ is an $n\times n$ matrix of continuous functions, $f:[0,\infty)\times D\to\R^n$ is continuous,  $D\subset C([-\tau,0];\R^n)$ is equipped with the uniform convergence metric, $\tau>0$ is the time-delay, and, as usual, $x_t$ denotes the past history of the system on the interval $[t-\tau,t]$, i.e., $x_t(\th)=x(t+\th)$ for $\th\in [-\tau,0]$.  The function $f$ is required to satisfy $f(t,0)=0$ for $t\ge 0$, and the ordinary differential equation (ODE)   $x'(t)=A(t)x(t)$  to be cooperative, thus  $a_{ij}(t)\ge 0$ must hold  for all $i\ne j$ and $t\ge 0$, where $a_{ij}(t)$ are  the entries of $A(t)$.
We restrict our setting to a
 class of delayed perturbations $f(t,x_t)$ with multiple time-varying discrete delays, having the particular form
\begin{equation}\label{i2}f(t,\phi)=(f_1(t,\phi_1),\dots, f_n(t,\phi_n))
\end{equation}
for $t\ge 0$ and  $\phi=(\phi_1,\dots,\phi_n)\in D$, where
$f_i(t,\phi_i)= \sum_{k=1}^m n_{ik}(t,\phi_i(-\tau_{ik}(t)))$ and $n_{ik}(t,x),\tau_{ik}(t)$ are  continuous, bounded and nonnegative  functions, for all $i,k$. For simplicity, this paper deals  with discrete delays only; however,  as pointed out later in Section 3, straightforward generalizations to some  families of perturbations with distributed delays are possible.

Inserting \eqref{i2} in \eqref{i1} leads to  systems of the form
\begin{equation}\label{i3}
x_i'(t)=\sum_{j=1}^n a_{ij}(t)x_j(t)+\sum_{k=1}^m   n_{ik}(t,x_i(t-\tau_{ik}(t))),\  i=1,\dots,n,\ t\ge 0,
\end{equation}
which can be interpreted as a structured population model for $n$ populations,
 see Section 3 for an additional set of hypotheses, as well as for some biological elements  of the model.

 \smal

 In the present paper, the main idea is to take full advantage of the properties of the cooperative  non-delayed linear system $x'(t)=A(t)x(t)$, to further analyse the large-time behaviour of solutions of system \eqref{i3}.
We shall impose conditions on the coefficients of  the linear system $x'(t)=A(t)x(t)$, in order to have its global exponential stability. This property and the monotonicity of $x'(t)=A(t)x(t)$ will play an important role in the study of  \eqref{i3}.
 Although the  nonlinearities  \eqref{i2} are in general non-monotone,    the techniques  exploited here are largely  based on  results  of comparison of solutions  (see \cite{Smith}), applied to some convenient  auxiliary cooperative DDE systems. This method is used to address  the global asymptotic behaviour of solutions of system \eqref{i3}, in what concerns its dissipativity,  uniform persistence and the global asymptotic stability of the null solution.
To some extent and in different frameworks, similar techniques have inspired the papers \cite{Faria14b, Faria15, Liu10,Liu11, Wang13}.  Some relevant applications are given. We also hope that the present  results can be used  to further address other aspects of the global dynamics of \eqref{i3}.

\smal

As a  significant example  of systems in the form \eqref{i3},
we shall  consider a
{\it non-autonomous} Nicholson system
with patch structure and multiple time-dependent discrete delays, given by
\begin{equation}\label{i4}
x_i'(t)=-d_i(t)x_i(t)+\sum_{j=1,j\ne i}^n a_{ij}(t)x_j(t)+\sum_{k=1}^m \be_{ik}(t)    x_i(t-\tau_{ik}(t))e^{-c_{ik}(t)x_i(t-\tau_{ik}(t))},\,  i=1,\dots,n,
\end{equation}
where all the coefficient and delay functions  are continuous, nonnegative and bounded. We stress that results concerning  multi-dimensional  versions of the famous Nicholson's blowflies equation \cite{GBN}  are still  quite limited, with most authors  treating only   {\it autonomous}  systems.

The papers of   Faria and R\"ost \cite{FariaRost}, on {\it autonomous} Nicholson systems, and of Obaya and Sanz  \cite{ObayaSanz}, on uniform and strict persistence for monotone skew-product semiflows, were
a strong motivation for the  present work. Here, the authors  further pursue  their previous research,   and extend it to general non-autonomous Nicholson systems: in fact, we aim to obtain results on  extinction, uniform persistence and  permanence of \eqref{i4} as simple  illustrations of our main results, proven for a much larger  family of DDEs of the form \eqref{i3}.

 \med

The contents  of the paper are now  briefly described. In Section 2, we study a cooperative ODE $x'(t)=A(t)x(t)$ and give sufficient conditions for its global exponential stability; in this case, a system of the form \eqref{i1} is  dissipative if the delayed perturbation $f(t,x_t)$ is bounded. In Section 3, we start by introducing  a set of assumptions for a family of  DDEs \eqref{i1}, give some  biological interpretation of the models and refer to some recent literature. The main results are then presented, providing very general criteria for both the global asymptotic stability of the trivial solution (in biological terms, this  implies the extinction of the populations in all patches) and the uniform persistence of such systems. A comparison with  results in  \cite{Faria15, Liu13,Wang13, Zhang}  is also given, and some questions are raised to be left as open problems.
Finally, in Section 4 we consider systems with autonomous coefficients (but with possible   time-dependent delays):  from the results in Section 3 and by a careful analysis of properties of cooperative matrices, we provide necessary and sufficient conditions for both their permanence and extinction, even in the case of a reducible community matrix. These sharp criteria improve and extend results for autonomous systems proven in recent literature. As an important example of application, throughout the paper our results are widely  illustrated with versions of  the Nicholson system \eqref{i4}.

\section{Preliminaries}

\setcounter{equation}{0}

In this section, we establish some preliminary results on stability for non-autonomous  linear homogeneous systems of ODEs of cooperative type. Although such systems have been widely studied (see e.g. \cite{Coppel,Fink,Hale69}), some optimal conditions for their asymptotic stability and global exponential stability are given here. For completeness of the reader, the authors opt to include  these conditions here, with the proof of a result whenever its precise statement   could not be found  elsewhere. We start with some standard definitions from the literature \cite{Coppel,Fink,Hale69}.

Consider an $n$-dimensional ODE   $x'=f(t,x)$  with  $f:[\al,\infty)\times D\subset \R^{n+1}\to \R^n$ sufficiently regular so that for any $(t_0,x_0)\in  [\al,\infty)\times D$ there exists a unique solution, denoted by $x(t,t_0,x_0)$, of the initial value problem $x'=f(t,x), x(t_0)=x_0$, defined on $[\al,\infty)$.  To simplify the writing, let $D=\R^n$. We further assume that $x=0$ is a solution, i.e., $f(t,0)=0, t\ge \al$. The zero  solution  is said to be {\it stable} on the interval $[\al,\infty)$ if for any $\vare>0$ and $t_0\ge \al$  there is $\de=\de(\vare,t_0)>0$ such that  $|x(t,t_0,x_0)|< \vare$ for all $t\ge t_0$, whenever $|x_0|< \de$; $x= 0$ is {\it uniformly stable} if it is stable and $\de$ above can be chosen independently of $t_0\ge \al$. The zero  solution is said to be {\it uniformly asymptotically stable} on $[\al,\infty)$ if it is uniformly stable  and there is $b>0$ such that, for any $\vare>0$, there is $T=T(\vare)>\al$ such that, for any $t_0\ge \al$ and $|x_0|< b$, we have  $|x(t,t_0,x_0)|< \vare$ for all $t\ge t_0+T$; and $x=0$ is  {\it globally exponentially  stable} on  $[\al,\infty)$ if there exist $K,\be>0$ such that $|x(t,t_0,x_0)|\le Ke^{-\be(t-t_0)}|x_0|$ for all $t\ge t_0\ge \al$ and $x_0\in\R^n$.  In general, the stability of a particular solution $\tilde x(t)$ of an ODE   $x'=f(t,x)$ is defined as the stability of the zero solution of  $y'=f(t,\tilde x+y)-f(t,\tilde x).$

The usual partial order in $\R^n$ relative to the cone $[0,\infty)^n$ is denoted here by $\le$, i.e.,  for $x,y\in\R^n$, $x\le y$ means $y-x\in [0,\infty)^n$; we write $x\ll y$ whenever $y-x\in (0,\infty)^n$. The notations $\ge$ and $\gg$ have then a clear meaning. In particular,
a vector $v$  in $\R^n$ is said to be  {\it positive} ({\it nonnegative})  if all its components are positive (nonnegative), and we write $v\gg 0$ ($v\ge 0$); by $v>0$ we mean that $v\ge 0$ and $v\ne 0$.

\begin{lem}\label{lem1}  Consider a non-autonomous linear ODE
\begin{equation}\label{01}x'(t)=A(t)x(t),\q t\ge \al,\end{equation}
where $\al\in\R$ and $A(t)=[a_{ij}(t)]$ is an $n\times n$ matrix of   functions such that:
\vskip 0mm
(a1) $a_{ij}$ are continuous on $[\al,\infty)$, $a_{ij}(t)\ge 0,i\ne j, a_{ii}(t)<0$ for all $t\ge \al$ and $i,j\in\{1,\dots,n\}$;
\vskip 0mm
(a2) there exists a vector $v=(v_1,\dots,v_n)\gg 0$ such that $A(t)v\le 0$ for all $t\ge \al$.
\vskip 0mm
Then, for any solution $x(t)$ of \eqref{01}, $|x(t)|_{v^{-1}}$ is non-increasing on $t\in[\al,\infty)$, where $|\cdot|_{v^{-1}}$ is the norm in $\R^n$ defined by $|x|_{v^{-1}}=\max_{1\le i\le n} (v_i^{-1}|x_i|)$ for $x=(x_1,\dots,x_n)$.
\end{lem}

\begin{proof}  Rescaling the variables by $\hat x_i(t)=v_i^{-1}x_i(t)\, (1\le i\le n)$, where $v=(v_1,\dots,v_n)\gg 0$ is a vector as in (a2), we obtain a new linear ODE $\hat x'(t)=\hat A(t) \hat x(t)$, where the matrix $\hat A(t)=[\hat a_{ij}(t)]$ has entries $\hat a_{ij}(t)=v_i^{-1}a_{ij}(t)v_j$. In this way, and after dropping the hats for simplicity, we may consider \eqref{01} where $v={\bf 1}:=(1,\dots,1)$ is the positive vector  in (a2) and
$|x|_{v^{-1}}=\max_{1\le i\le n} |x_i|$.

Let $x(t)\ne 0$ be a solution  of \eqref{01}. To prove the claim, we show that $|x(t)|$ is non-increasing on each fixed interval $J=[t_0,t_1],\ \al\le t_0<t_1$. Define $u_j=\max_J |x_j(t)|$, and let $u_i=\max_{1\le j\le n}u_j$, with $u_i=|x_i(t_*)|$ for some $t_*\in J$. It is sufficient to show that $u_i=|x_i(t_0)|$.

We suppose that $x_i(t_*)>0$; the case $x_i(t_*)<0$ is treated in a similar way. Denoting $d_i(t)=-a_{ii}(t)$ and $D_i(t)=\int_{t_0}^t d_i(s)\, ds$, for $t\in J$ we have
$x_i'(t)+d_i(t)x_i(t)\le d_i(t)u_i$. Hence
$$x_i(t)\le x_i(t_0)e^{-D_i(t)}+u_i(1-e^{-D_i(t)}),\q t\in J.$$
In particular for $t=t_*$ we derive $u_ie^{-D_i(t_*)}\le x_i(t_0)e^{-D_i(t_*)}$, thus $u_i=x_i(t_0)$.
\end{proof}

 \begin{lem}\label{lem2}  For the linear ODE system \eqref{01},  assume
\vskip0mm
(a1') $a_{ij}$ are  uniformly  continuous and bounded on $[\al,\infty)$, $a_{ij}(t)\ge 0,i\ne j, a_{ii}(t)<0$ for all $t\ge \al$ and $i,j\in\{1,\dots,n\}$;
\vskip 0mm
(a2') there exists a vector $v=(v_1,\dots,v_n)\gg 0$ such that $A(t)v\le 0$ for all $t\ge \al$, and $\liminf_{t\to\infty}A(t)v\ll 0$, in the sense that there exists a sequence $ t_k\to \infty$ such that $\lim_k(A(t_k)v)_i<0,\, i=1,\dots,n$.
\vskip 0mm
Then,  \eqref{01} is  asymptotically stable; in other words, $x=0$ is stable and $\lim_{t\to\infty}x(t)=0$, for all solutions of \eqref{01}.
\end{lem}

\begin{proof} As in the above proof and without loss of generality, consider $v={\bf 1}$  in (a2') and the norm
$|x|=\max_{1\le i\le n} |x_i|$ in $\R^n$.
From Lemma \ref{lem1}, \eqref{01} is uniformly stable. We now prove that  the trivial solution  is a global attractor of all solutions.

Let $x(t)\ne 0$ be a solution of \eqref{01}, and define $c=\lim_{t\to\infty}|x(t)|$. We want to show that $c=0$.
 In order to obtain a contradiction, suppose that $c>0$. By (a2'), take $t_k\to\infty$ such that $\al_i:=\lim_k(-d_i(t_k)+\sum_{j\ne i}a_{ij}(t_k))<0$ for all $i$. In particular for such a sequence, $|x(t_k)|\searrow c$, and thus there exists $i\in \{1,\dots,n\}$ and a subsequence, still denoted by $(t_k)$, such that either $x_i(t_k)=|x(t_k)|\to c$ or $x_i(t_k)=-|x(t_k)|\to -c$. We  only consider the situation $x_i(t_k)\to c$ for some $i$, the other is treated in a similar way.
 We now consider separately two cases.

  First, suppose that there exists $\lim_{t\to \infty} x_i(t)=c.$
 Since the entries of $A(t)$ are bounded and uniformly continuous and $x(t)$ is uniformly bounded on $[\al,\infty)$, one easily shows that all components  $x_j(t)$ and $x_j'(t)$ are uniformly continuous  on $[\al,\infty)$. From the Barbalat Lemma, we derive that there is $\lim_{t\to\infty} x_i'(t)=0$, and in particular obtain  $\lim_k x_i'(t_k)=0$.
On the other hand, from \eqref{01} we have
 $$x_i'(t_k)\le  x_i(t_k)\Big [-d_i(t_k)+\sum_{j\ne i}a_{ij}(t_k)\Big].$$
 Taking  limits, the above inequality leads to $0\le c\al_i$, which is a contradiction.

 Next, consider the case when  $\ul x_i:=\liminf_{t\to \infty} x_i(t)<\limsup_{t\to \infty} x_i(t)=c.$
  From the inequality above, it is clear that $t_k$ are not local extrema points, since $x_i'(t_k)<0$.  However, by reducing to a subsequence if necessary, we may consider that each $t_k$ lies  between a local maximum point $\ol t_k$ and a local minimum point $\ul t_k$,  to its left and to its right respectively, with $x_i(\ol t_k)\to c, x_i(\ul t_k)\to \ul x_i$. We may take a sequence $(s_k)$ such that, for all $k\in\N$,
  $$t_k<\ul t_k<s_k<\ol t_{k+1}<t_{k+1},\q x_i(s_k)=x_i(t_k),\q x_i'(s_k)\ge 0.$$
  Since the map $t\to |x(t)|$ is non-increasing, we have $x_j(s_k)\le |x(s_k)|\le |x(t_k)|=x_i(t_k),\, 1\le j\le n.$
We derive
  $$0\le x_i'(s_k)\le x_i(t_k)\Big [-d_i(s_k)+\sum_{j\ne i}a_{ij}(s_k)\Big]\le 0,$$
 and therefore $ x_i(t_k)=0$, which is not possible.
  \end{proof}

When $f(t,x)$ is periodic in $t$, a solution of $x'=f(t,x)$ is uniformly asymptotically stable if it is asymptotically stable. This is not true if periodic is replaced by almost periodic (see \cite{Fink},~p.~191 for a counter-example). 
Moreover,  for  a linear system $x'=A(t)x$, where $A(t)$ is an $n\times n$ matrix of continuous functions, it is well known that the concepts of  global  exponential  stability and uniform asymptotic stability on an interval $[\al,\infty)$ are equivalent (see \cite{Coppel,Hale69}). 
Therefore, the  following criterion  is straightforward for periodic systems, however it applies to the more general case of almost periodic linear systems.

  \begin{thm}\label{thm20} Let $A(t)=[a_{ij}(t)]$ be an $n\times n$ matrix  of  almost periodic functions on $\R$ satisfying (a1), (a2) on $\R$, with $A(t_0)v\ll 0$ for some $t_0\in \R$.
  Then, \eqref{01} is globally exponentially stable.\end{thm}

\begin{proof}
Let $H(A)$ be the hull of $A$, that is, the closure for the topology of uniform convergence of the set of shifted maps $\{\theta_t A({\cdot})=A({\cdot}+t)\mid t\in \R\}$  \cite{Coppel,Fink}. $H(A)$ is   a compact metric space.
Since $A$ is almost periodic, it follows that $A$ satisfies (a1'), (a2').  The orbit $\{\theta_t A\mid t\in  \R\}$ is dense in the hull, thus actually any $B\in H(A)$ satisfies (a1'), (a2') as well. By Lemma \ref{lem2}, all solutions of all the systems $x'=B(t)x$, with $B\in H(A)$, tend to $0$ as $t\to \infty$. At this point, the spectral theory of Sacker and Sell \cite{sase} applies and permits  to conclude that \eqref{01} is globally exponentially stable.
 \end{proof}



Usually, the global exponential  stability of \eqref{01} is  obtained by assuming that $A(t)$ is {\it strongly} uniformly row (or column) dominant. The theorem below follows from Proposition 6.3 in \cite{Coppel}.

\begin{thm}\label{thm21} Consider an $n\times n$ matrix $A(t)=[a_{ij}(t)]$ of  bounded continuous functions satisfying (a1),  and suppose that \vskip 0mm
(a3) there exist a vector $v=(v_1,\dots,v_n)\gg 0$ and $T\ge \al,\de>0$ such that $(A(t)v)_i\le -\de$ for all $t\ge T,\, i=1,\dots,n$.
\vskip 0mm
  Then, \eqref{01} is globally exponentially  stable.\end{thm}

%


   \begin{rmk}\label{rmk2.1}  {\rm For any fixed $t$, the matrix $-A(t)$ is a non-singular M-matrix if and only if there exists a positive vector $v=v(t)$ such that  $A(t)v\ll 0$;  thus, condition (a3) above not only demands that  $-A(t)$ are non-singular M-matrices, for $t$ sufficiently large, but also that there exist positive vectors $v,\eta$, which do not depend on $t$, such that $A(t)v\le -\eta$ (see \cite{Fiedler} and Section 4 for more details on M-matrices).}
  \end{rmk}

For $\tau\ge 0$, consider the Banach space $C:=C([-\tau,0];\R^n)$  equipped with the norm
$\|\phi\|=\max_{\th\in[-\tau,0]}|\phi(\th)|$, where $|\cdot|$ is  a fixed norm in $\R^n$.  The case of no delays ($\tau=0$) is included, in which case $C$ is identified with $\R^n$. We now consider    DDEs obtained by adding a  bounded delayed perturbation   $f(t,x_t)$ to systems \eqref{01}, where, as before,  $x_t\in C$ is given by $x_t(\th)=x(t+\th), -\tau\le \th\le 0$.
 For  simplicity, in what follows we take $\al=0$, but any $\al\in\R$ could be considered.

From Theorem \ref{thm21}, one obtains:

\begin{thm}\label{thm1} Consider an $n\times n$ matrix $A(t)=[a_{ij}(t)]$ of bounded functions satisfying (a1), (a3) on $[0,\infty)$, and a function $f:[0,\infty)\times C\to \R^n$  continuous and bounded. Then, all solutions of the DDE
\begin{equation}\label{02}
x'(t)=A(t)x(t)+f(t,x_t),\q t\ge 0,
\end{equation}
are defined  on $[0,\infty)$ and \eqref{02}
is dissipative, i.e., there exists $M>0$ such that $\dps \limsup_{t\to\infty}|x(t)|\le M$ for any solution $x(t)$ of \eqref{02}.
\end{thm}

\begin{proof}  Let $|f(t,\var)|\le L$ for $t\ge 0,\var \in C$. From Theorem \ref{thm21}, there are $K>0,\al>0$ such that
$|X(t)X^{-1}(t_0)|\le Ke^{-\al (t-t_0)},\, t\ge t_0\ge 0,$
where $X(t)$ is a fundamental  solution matrix for \eqref{01}. By the variation of constants formula, the solutions $x(t)$ of \eqref{02} satisfy
\begin{equation}\label{2.3}
x(t)=X(t)X^{-1}(t_0)x(t_0)+X(t)\Big (\int_{t_0}^t X^{-1}(s) f(s,x_s)\, ds\Big)\q (t,t_0\ge 0),
\end{equation}so that
$|x(t)|\le Ke^{-\al (t-t_0)} |x(t_0)|+\frac{KL}{\al}(1-e^{-\al(t-t_0)})\to \frac{KL}{\al}\q {\rm as}\q t\to\infty.$
\end{proof}

We now set some further notation.  Let $C^+$ be the cone of nonnegative functions in $C$, $C^+=C([-\tau,0];[0,\infty)^n)$, and $int \, C^+$ its interior. Hereafter, $\le $  also denotes the usual partial order generated by  $C^+$: $\phi\le \psi$ if and only if $\psi-\phi\in C^+$; by $ \phi\ll \psi $,   we mean that $\psi-\phi\in int\, C^+$. The definition of the relations $\ge$ and $\gg $ are then clear; thus,  we write $\psi\ge 0$ for $\psi\in C^+$ and $\psi\gg 0$ for $\psi \in int \, C^+$.  
A vector $v\in \R^n$  is identified in $C$ with the constant function $\psi(s)=v$ for $-\tau\le s\le 0$.

Let $D\subset C([-\tau,0];\R^n)\ (\tau\ge 0)$ be open, and consider a  non-autonomous DDE written as
\begin{equation}\label{2.4}x'(t)=f(t,x_t),\q t\ge 0,\end{equation}
where $f:[0,\infty)\times D\to \R^n$ is continuous and regular enough  so that
the initial value problem is well-posed, in the sense that for each $(\sigma,\phi)\in [0,\infty)\times D$ there exists a unique solution of  the problem $x'(t)=f(t,x_t), x_\sigma=\phi$,  defined on a maximal interval of existence. This solution will be denoted by $x(t,\sigma,\phi)$ in $\R^n$ or  $x_t(\sigma,\phi)$ in $C$.   When considering more than one DDE $x'(t)=f(t,x_t)$, the notation $x(t,\sigma,\phi,f)$  where the argument $f$ is made explicit  will be  used  to clarify  which DDE is being considered.

To simplify the terminology, we say that  \eqref{2.4} is {\it cooperative} if it satisfies Smith's {\it quasi-monotone condition} (Q), given by (see \cite{Smith})

\vskip 3mm

{(Q)} for $\phi,\psi\in D,\phi\le \psi$ and $\phi_i(0)=\psi_i(0)$, then $f_i(t,\phi)\le f_i(t,\psi), \ i=1,\dots,n,t\ge 0$.

\vskip 3mm


It is well-known that (Q) guarantees monotonocity of solutions relative to initial data and allows   comparison of solutions between two related DDEs,  $x'(t)=f(t,x_t), x'(t)=g(t,x_t)$ with $f\le g$:   if at least one of them is cooperative, then
$x(t,\sigma, \phi,f)\le x(t,\sigma, \psi,g)$ for $t\ge \sigma$ if $\phi\le \psi$ (\cite{Smith}).
These and other properties of cooperative ODEs and DDEs
  will turn out to be very useful in the next sections.
 The lemma below will  be often applied, see p.~82 of~\cite{Smith}.

\begin{lem}\label{lem3} \cite{Smith} Consider \eqref{2.4}  in $D\subset C([-\tau,0];\R^n)$, and let $v=(v_1,\dots,v_n)\in \R^n$.\vskip 0cm
(i) If  $f_i(t,\phi)\le 0$ for all $ i=1,\dots,n,t\ge 0$ whenever $\phi\in D,\phi\le v$ and $\phi_i(0)=v_i$, then the set $\{\phi\in D: \phi\le v\}$ is positively invariant for \eqref{2.4}.
\vskip 0cm
(ii) If  $f_i(t,\phi)\ge 0$ for all $i=1,\dots,n,t\ge 0$ whenever $\phi\in D,\phi\ge v$ and $\phi_i(0)=v_i$, then the set $\{\phi\in D: \phi\ge v\}$ is positively invariant for \eqref{2.4}.
\end{lem}

\begin{rmk}\label{rmk2.2} {\rm Clearly, if (a1)
 is satisfied,  then \eqref{01} is a cooperative system and the nonnegative cone $[0,\infty)^n$ is forward invariant. If in addition (a2) is satisfied and $v=(v_1,\dots,v_n)\gg 0$ is as in (a2), for $x\in\R^n$  such that $x\le v$ and $x_i=v_i$, then $(A(t)x)_i\le 0$. This  implies that the interval $[0,v]:=[0,v_1]\times \cdots\times[0,v_n]$ is forward invariant as well.}
\end{rmk}

\section{Global behaviour for a class of  non-monotone and non-autonomous DDEs}
\setcounter{equation}{0}

In this section, we consider $n$-dimensional delayed  structured models \eqref{i1}, where the
  linear ODE system \eqref{01}  is globally exponentially  stable, $f$ is  continuous,   bounded, and, in general, non-monotone. Although some generalizations are possible, we restrict our framework to  perturbations   $f(t,x_{t})=(f_1(t,x_{1,t}),\dots, f_n(t,x_{n,t}))$,
 with each component $f_i(t,\phi_i)$  of the form
$f_i(t,\phi_i)= \sum_{k=1}^m n_{ik}(t,\phi_i(-\tau_{ik}(t))$, for $t\ge 0, \phi=(\phi_1,\dots,\phi_n)\in C$.
Moreover, we suppose that  $n_{ik}(t,0)=0$ for $t\ge 0$ and have partial derivative with respect to the second variable at $x=0^+$ given by $\frac{\p n_{ik}}{\p x}(t,0)=\be_{ik}(t)\ge 0$; thus  $n_{ik}(t,x)$ is written as $n_{ik}(t,x)=\be_{ik}(t)h_{ik}(t,x)$ with $h_{ik}(t,0)=0,\frac{\p h_{ik}}{\p x}(t,0)=1,t\ge 0$.
Below, some additional assumptions  on $n_{ik}(t,x)$ will be imposed.
This leads to a non-autonomous  system with 
multiple discrete  time-dependent delays of the form
\begin{equation}\label{3}
x_i'(t)=-d_i(t)x_i(t)+\sum_{j=1,j\ne i}^n a_{ij}(t)x_j(t)+\sum_{k=1}^m \be_{ik}(t)    h_{ik}(t,x_i(t-\tau_{ik}(t))),\  i=1,\dots,n,\ t\ge 0.
\end{equation}

 Throughout the remainder of this paper,   either the whole or a part of the following set of hypotheses  will be imposed:
\begin{itemize}
\item[(h1)] the  functions $d_i,a_{ij}\, (j\ne i)$ are   continuous and bounded, with $a_{ij}(t)\ge 0,i\ne j, d_i(t)>0$ for $t\ge 0$ and $i,j\in\{1,\dots,n\}$
\item[(h2)] there exist a vector $v=(v_1,\dots,v_n)\gg 0$ and $\de >0,T_0\ge 0$   such that
$d_i(t)v_i\ge \sum_{j=1,j\ne i}^n a_{ij}(t)v_j+\de $ for $t\ge T_0, i\in\{1,\dots,n\}$;


\item[(h3)]  $\tau_{ik},\be_{ik}$ are continuous and bounded, with $\tau_{ik}(t)\ge 0,\be_{ik} (t)\ge 0$ and
$$\be_i(t):=\sum_{k=1}^m\be_{ik} (t)>0$$ for $t\in [0,\infty)$,  $i\in \{1,\dots,n\},k\in \{1,\dots,m\}$;

\item[(h4)] $h_{ik}:[0,\infty)\times [0,\infty)\to [0,\infty)$ are  bounded, continuous, $h_{ik}(t,x)$ are  locally Lipschitzian in $x$, with
$$h_i^-(x)\le h_{ik}(t,x)\le h_i^+(x),\q {\rm}\q t,x\ge 0,k=1,\dots,m,$$ where $h_i^\pm:[0,\infty)\to [0,\infty)$ are  continuous on $[0,\infty)$ and continuously differentiable in a vicinity of $0^+$,  with $h_i^\pm(0)=0,(h_i^\pm)'(0)=1$ and $h_i^-(x)>0$ for $x>0$,  $i\in \{1,\dots,n\}$.
\end{itemize}

For simplicity, here we only treat
 non-autonomous  systems with discrete non-autonomous delays, but our framework applies with straightforward   adjustments to
 the more general case of systems with multiple distributed  time-varying delays of the form
\begin{equation}\label{003}
x_i'(t)=-d_i(t)x_i(t)+\sum_{j=1,j\ne i}^n a_{ij}(t)x_j(t)+  f_i(t,x_{i,t}), \  i=1,\dots,n,
\end{equation}
with  $f_i(t,x_{i,t})$ given by
\begin{equation}\label{0003}
f_i(t,x_{i,t})=\sum_{k=1}^m \be_{ik}(t) h_{ik}\big (t, L_{ik} (t,x_{i,t})\big )\ {\rm or}\
f_i(t,x_{i,t})=\sum_{k=1}^m \be_{ik}(t)  L_{ik}\big (t, h_{ik}(\cdot, x_{i,t}(\cdot))\big ),
\end{equation}
where
$$L_{ik}(t,\phi)=\int_{-\tau}^0 \phi(s)\,  d_s\eta_{ik}(t,s)\q {\rm for}\  t\ge 0,\, \phi\in C([-\tau,0],\R),$$
 $\tau>0$,  the measurable functions $\eta_{ik}:[0,\infty)\times[-\tau,0]\to \R$ are continuous on $t$, with $\eta_{ik}(t,\cdot)$ non-decreasing and normalized so that
$\int_{-\tau} ^0 d_s\eta_{ik}(t,s)= 1 ,\ i=1,\dots,n,k=1,\dots,m,t\ge 0,$
and  for which (h3), (h4) hold. Besides  \eqref{0003}, and under some natural conditions, other forms of dependence on distributed delays can  be incorporated in \eqref{003}.

\vskip 1mm

In what follows, we   refer to the $n\times n$ matrix-valued functions defined on $[0,\infty)$  by
\begin{equation}\label{3.2}
\begin{split}
D(t)=diag\, (d_1(t),\dots,d_n(t)),&\q  A(t)=[a_{ij}(t)]\\
 B(t)=diag\, (\be_1(t),\dots,\be_n(t)),&\q
 M(t)=B(t)+A(t)-D(t),\qq t\ge 0,
\end{split}
\end{equation}
where $a_{ii}(t)\equiv0$. The matrix $M(t)$ is often designated as the {\it community matrix} of the population system \eqref{3}.

\begin{rmk}\label{rmk3.0} {\rm
We stress that under (h1), (h2) the linear homogeneous ODE
$x'(t)=-[D(t)-A(t)]x(t)$
possesses  two important features: it is cooperative and globally exponentially  stable. Of course, if $A(t)$ is periodic or almost periodic, Theorem \ref{thm20} allows us to replace (h2) by the weaker condition $[D(t)-A(t)]v\ge 0$ for $t\in \R$ and $[D(t_0)-A(t_0)]v\gg 0$ for some  $t_0\in \R$ and   $v\gg 0$.}
\end{rmk}

 System \eqref{3} can be interpreted as a  model  for $n$  populations  structured into $n$ classes or patches, with migration among them: $x_i(t)$ denotes the density of the $i$th population; $a_{ij}(t)$ is the migration rate of the population in class $j$ moving to class $i$; $d_i(t)$ is the coefficient of instantaneous loss for class $i$, which incorporates both the death rate and the emigration rates of the population  that leaves class $i$ to move to other classes; the birth contribution for each population is given by the  nonlinear terms $\sum_k \be_{ik}(t)h_{ik}(t, x_i(t-\tau_{ik}(t)))$.

 With this interpretation, $d_i(t)=m_i(t)+\sum_{j\ne i} a_{ji}(t)$, where
$m_i(t)$ is the death rate for the $i$th population, so  it is natural to
impose $a_{ij}(t)\ge 0$ and   $d_i(t)>\sum_{j\ne i} a_{ji}(t)$ for all
$i,j$, i.e., $D(t)-A(t)^T$ is uniformly diagonally dominant for $t\ge 0$. It is also natural to assume that $a_{ij}(t)=\vare_{ij}(t)a_{ji}(t)$
for $i\ne j$ and $t\ge 0$, with $\vare_{ij}(t)\in (0,1]$, to account for some loss of  the  populations, when moving to  different patches (see \cite{Takeuchi}), thus
$[D(t)-A(t)]{\bf 1}\gg 0$ for $t\ge 0$.  If the mortality rates $m_i(t)$ are bounded below by a positive constant $m_0$, then $[D(t)-A(t)]{\bf 1}\ge m_0{\bf 1}$ for $t\ge 0$.
To some degree, these comments justify assumption  (h2) from a biological  point of view.

  Following the general approach in the literature,  here multiple (time-varying) discrete delays have been introduced in the birth function.
  In biological terms, most situations do not require the consideration of more than one delay, either a discrete or a distributed delay, but occasionally multiple delays should be incorporated  in each equation. For  examples of such situations, we refer   to generalizations of the classic Mackey-Glass model for the production of red blood cells in \cite{Belair} and to \cite{Smith11} for other   references.

\vskip 1mm

As an important example of application, we have in mind the following  non-autonomous Nicholson system
with patch structure and multiple time-dependent discrete delays:
\begin{equation}\label{04}
x_i'(t)=-d_i(t)x_i(t)+\sum_{j=1,j\ne i}^n a_{ij}(t)x_j(t)+\sum_{k=1}^m \be_{ik}(t)    x_i(t-\tau_{ik}(t))e^{-c_{ik}(t)x_i(t-\tau_{ik}(t))},
\end{equation}
for $ i=1,\dots,n,\ t\ge 0$. For \eqref{04}, we shall always assume that the coefficient and delay functions satisfy (h1), (h3) and that $c_{ik}(t)\ge c_i>0$ are continuous and bounded.  With nonlinearities given by $h_{ik}(t,x)=xe^{-c_{ik}(t)x}$ for all $i,k$,   (h4) is obviously satisfied.

The {\it autonomous} version of \eqref{04} with $n=1$ and $m=1$ is  the famous  {\it Nicholson's blowfly equation}, given by
$N'(t)=-dN(t)+\be N(t-\tau)e^{-aN(t-\tau)}\ (d,\be, a,\tau>0)$. A large-scale literature on the scalar Nicholson's blowflies equation, on a number of generalizations and on related models has been produced since its introduction  by Gurney et al.~\cite{GBN},   and real world applications  implemented. Nevertheless, a number of   problems regarding scalar Nicholson-type equation  still remain unsolved, see \cite{BerBrav,BBI} and references therein.
On the other hand, results concerning  multi-dimensional  versions of such models  are still  quite limited.
Not only is the literature on Nicholson systems very sparse, but also most authors have only treated   {\it autonomous} Nicholson systems, and only recently have non-autonomous Nicholson systems  been considered. See \cite{BIT, Faria11, FariaRost, Liu09, Liu10, Liu11, Wang13, WWC,Zhou}, also for  biological details of the models and additional references. 

Besides  Ricker-type nonlinearities as  in the non-autonomous Nicholson system  \eqref{04}, other useful population models can be written in  the form \eqref{3}. Among them, are models with 
Mackey-Glass type nonlinearities of the form (see \cite{MG})
$$h_{ik}(t,x)=xe^{-c_{ik}(t)x^\al}\q (\al>0)\q {\rm or}\q
h_{ik}(t,x)=\frac{x}{1+c_{ik}(t)x^\al}\q (\al\ge 1),$$
which satisfy (h4) if $c_{ik}(t)$ are continuous and bounded below and above by positive constants.

\vskip 1mm

System \eqref{3} is considered as a DDE in   $C=C([-\tau,0];\R^n)$, where $\tau=\max_{i,k}\sup_{t\ge 0} \tau_{ik}(t)$.  Unless specifically mentioned,
$\|\phi\|=\max_{\th\in[-\tau,0]}|\phi(\th)|$ for $\phi\in C$, where $|\cdot|$ is the maximum norm in $\R^n$.
Motivated by the applications to mathematical biology,
only nonnegative solutions of \eqref{3} are meaningful. For this reason, initial conditions are taken in either $C^+$ or $C_0$, where $$C_0=\{ \var\in C^+: \var(0)\gg 0\}.$$

Together with \eqref{3}, we also consider its linearization  at the origin:
\begin{equation}\label{7}
y_i'(t)=-d_i(t)y_i(t)+\sum_{j=1,j\ne i}^n a_{ij}(t)y_j(t)+\sum_{k=1}^m \be_{ik}(t)    y_i(t-\tau_{ik}(t)),\q i=1,\dots,n.\end{equation}
Write \eqref{3}, \eqref{7} as $x'(t)=f(t,x_t), y'(t)=g(t,y_t)$ respectively, where $f=(f_1,\dots,f_n),g=(g_1,\dots,g_n)$ and
$$f_i(t,\phi)=-d_i(t)\phi_i(0)+\sum_{j=1,j\ne i}^n a_{ij}(t)\phi_j(0)+\sum_{k=1}^m \be_{ik}(t)   h_{ik}( t, \phi_i(-\tau_{ik}(t)))$$
and
$$g_i(t,\phi)=-d_i(t)\phi_i(0)+\sum_{j=1,j\ne i}^n a_{ij}(t)\phi_j(0)+\sum_{k=1}^m \be_{ik}(t)    \phi_i(-\tau_{ik}(t)).$$

 Assume (h1), (h3), (h4).
  For $t\ge 0$ and $\phi\ge 0,\phi_i(0)=0$, then $f_i(t,\phi)\ge 0$ and $g_i(t,\phi)\ge 0$, which implies that $x(t):=x(t,0,\phi,f)\ge 0$ and $y(t):=y(t,0,\phi,g)\ge 0$ for $t\ge 0$. Moreover, $x_i'(t)\ge -d_i(t)x_i(t)$ and $y_i'(t)\ge -d_i(t)y_i(t)$ for $t\ge 0$ and $1\le i\le n$. Hence both $C^+$ and $C_0$ are positively invariant for \eqref{3} and  \eqref{7}.
The next result  is a consequence of Theorem \ref{thm1}.

\begin{thm}\label{thm2} Under the assumptions (h1)-(h4), all solutions of \eqref{3} with initial conditions in $C_0$ are defined and  strictly positive on $[0,\infty)$; moreover,  there exists $L>0$ such that, for any $\phi\in C_0$, there is $T=T(\phi)>0$ such that
\begin{equation}\label{3.5}
0<x_i(t,0,\phi)< L\q {\rm for}\q t\ge T,\ i=1,\dots,n.
\end{equation}
\end{thm}

%


We  now introduce a notation often used for DDEs (cf.  \cite{Smith}, p.~82): if there is no possibility of misinterpretation with intervals of $\R$ or $\R^n$, for $v\in \R^n$ we also denote $[0,v]$ and $[v,\infty)$ as the subsets of $C$ given by $[0,v]=\{ v\in C^+:\var \le v\}$ and $[v,\infty)=\{ v\in C:\var \ge v\}.$

\begin{lem}\label{lem4}  Under (h1), (h3), system  \eqref{7} is cooperative, and the following holds:
\vskip 0cm
(i) If there exist a vector $v=(v_1,\dots,v_n)\gg 0$ and $T_0\ge 0$ such that $M(t)v\le 0$ for $t\ge T_0$, then the sets $[0,cv]\cap C_0$ (where $c>0$) are   invariant for \eqref{7} with $t\ge T_0$; in particular, the  solutions of \eqref{7} are uniformly stable.
\vskip 0cm
(ii) If there exist a vector $v=(v_1,\dots,v_n)\gg 0$ and $T_0\ge 0$ such that $M(t)v\ge 0$ for $t\ge T_0$, then the sets $[cv,\infty)\cap C_0$ (where $c>0$) are   invariant for \eqref{7} with $t\ge T_0$.
\end{lem}

\begin{proof} Since $g$ satisfies (Q), \eqref{7} is cooperative.
Let $M(t)v\le 0$ for $t\ge T_0$, for some strictly positive vector $v=(v_1,\dots,v_n)\in\R^n$ and some $T_0\ge 0$.  For $\phi\in C^+$ with $\phi\le v$, if $\phi_i(0)=v_i$ for some $i$, then
$g_i(t,\phi)\le \big( M(t)v\big )_i\le 0$ for $t\ge T_0$, proving that $[0,v]\cap C_0$ is  positively invariant (see Lemma \ref{lem3}); since the system is linear, for any positive constant $c$ the set $[0,cv]\cap C_0$ is  positively invariant as well. From the monotonicity, it follows that the solution $y=0$ of \eqref{7} is uniformly stable. The proof of (ii) is similar.
\end{proof}

\begin{defn} The trivial solution $x\equiv 0$   of \eqref{3}  is said to be {\bf stable}
if for any $\vare> 0$  there is $\delta=\de (\vare)>0$ such that $\|x_t(0,\phi)\|<\vare$ for all $\phi \in C_0$ with $\|\phi\|<\delta$ and $t\ge 0$;  $0$ is said to be
{\bf globally attractive} (in $C_0$)  if $x(t,0,\phi)\to 0$ as $t\to\infty$, for all  solutions of \eqref{3} with initial conditions $x_0=\phi \in C_0$;  $0$  is {\bf globally asymptotically stable (GAS)}  if it is stable and globally attractive.
\end{defn}

The next result gives sufficient conditions for  the  stability and global attractivity of the trivial equilibrium. When \eqref{3} refers to a population model, the global attractivity of 0 means the extinction of the populations in all patches.

\begin{thm}\label{thm3}  Assume (h1), (h3) and (h4) with  $0<h_i^+(x)< x,\, x>0,1\le i\le n$. Further suppose that:\vskip 0cm
(i) there exist  $v=(v_1,\dots,v_n)\gg 0$ and $T_0\ge 0$ such that $M(t)v\le 0$ for $t\ge T_0$;\vskip 0cm
(ii)  either $\liminf_{t\to\infty} \be_i(t)>0$ or  $\limsup_{t\to\infty} (M(t)v)_i<0$, for all $ i=1,\dots,n$.\vskip 0cm
\noindent{Then the trivial solution of \eqref{3} is GAS in $C_0$. }
\end{thm}

\begin{proof} From (ii), for each $i=1,\dots,n$, either $\be_i(t)\ge \ul \be_i>0$ for $t$ large or $(M(t)v)_i\le -\la_i<0$ for $t$ large.  In particular,  together with (h3), conditions (i) and (ii) imply (h2).

  For $\phi \in C_0, t_0\ge 0$ and $i\in\{1,\dots,n\}$,
it holds $f_i(t,\phi)\le g_i(t,\phi)$. In this way,  the solutions of \eqref{3} and \eqref{7} satisfy $x(t,t_0,\phi,f)\le y(t,t_0,\phi,g),t\ge t_0$. From  Lemma \ref{lem4}, the zero solution of \eqref{3} is stable. Now, we show that it attracts all solutions  with initial conditions in $C_0$.

 With $\hat x_j(t)=x_j(t)/v_j$, system \eqref{3} reads as
\begin{equation}\label{3'}
\hat x_i'(t)=- d_i(t)\hat x_i(t)+\sum_{j=1,j\ne i}^n \hat a_{ij}(t)\hat x_j(t)+\sum_{k=1}^m  \be_{ik}(t)    \hat h_{ik}(t,\hat x_i(t-\tau_{ik}(t))),\  i=1,\dots,n,\ t\ge 0,
\end{equation}
where $ \hat a_{ij}(t)= v_i^{-1}a_{ij}(t) v_j, j\ne i$, and $\hat h_{ik}(t,x)=v_i^{-1}h_{ik}(t,v_ix)$ satisfy (h4). Hence, without loss of generality we consider the original system \eqref{3} and  take $v={\bf 1}$  in (i),~(ii).

The solutions $x(t)=x(t,t_0,\phi,f)$  are bounded, so define $u_j=\limsup_{t\to \infty} x_j(t)$ and let $u_i=\max_{1\le j\le n} u_j$. If $u_i>0$, by the fluctuation lemma take a sequence $(t_k)$ with $t_k\to \infty$, $x_i(t_k)\to u_i, x'_i(t_k)\to 0$. For any small $\vare >0$ with $u_i-\vare>0$, for $k$  large we get $t_k\ge T_0+\tau$, $u_j-\vare\le x_j(t)\le u_j+\vare$ and $h_{ik}(t,x_i(t))\le \max_{x\in [0,u_i+\vare]}h_i^+(x)$, for $t\in [t_k-\tau, t_k]$. Thus,
\begin{equation*}
\begin{split}
 x_i'(t_k)&\le -d_i(t_k)(u_i-\vare)+\sum_{j\ne i} a_{ij}(t_k) (u_j+\vare)+\sum_p \be_{ip}(t_k)    h_i^+(x_i(t_k-\tau_{ip}(t_k)))\\
&\le-u_i\Big ( d_i(t_k)-\sum_{j\ne i} a_{ij}(t_k)\Big)+\be_{i}(t_k)     \max_{x\in [0,u_i+\vare]}h_i^+(x)+O(\vare).
\end{split}
\end{equation*}
Taking limits $k\to\infty, \vare\to 0^+$, we derive that
$$0\le \limsup_{t\to\infty} u_i\Big (\be_{i}(t)-d_i(t)+\sum_{j\ne i} a_{ij}(t)\Big)+ (  \max_{x\in [0,u_i]}h_i^+(x)-u_i ) \liminf_{t\to\infty} \be_i(t).$$ Since $\max_{x\in [0,u_i]}h_i^+(x)<u_i$ and one of the conditions in (ii) is satisfied, this
 is not possible. Therefore $u_i=0$, and the proof is complete.
\end{proof}

 For the definitions of persistence and permanence given below, see e.g.~\cite{SmithThi}.

\begin{defn}\label{def3.1} A set $S\subset C^+$ is  an {\it admissible set  of initial conditions} for  $x'(t)=f(t,x_t)$ if  any  solution $x(t,0,\phi)$ with initial condition $x_0=\phi\in S$  satisfies $x_t\in S$ for $t\ge 0$, whenever it is defined.
A DDE $x'(t)=f(t,x_t)$
 is said to be {\bf persistent}  in $S$, for $S$ an admissible set  of initial conditions,  if all solutions $x(t,0,\phi)$ with $\phi\in S$ are defined and bounded below away from zero on $[0,\infty)$, i.e.,
$\liminf_{t\to\infty}x_i(t,0,\phi)>0$ for all $ 1\le i\le n,\phi\in S;$
and $x'(t)=f(t,x_t)$ is {\bf uniformly persistent}  in $S$ if there is  $m>0$ such that $\liminf_{t\to\infty}x_i(t,0,\phi)\ge m$ for all $1\le i\le n,\phi\in S$.
The system is said to be {\bf permanent} in $S$  if it is dissipative and uniformly persistent; in other words, all solutions $x(t,0,\phi),\phi\in S$, are defined on $[0,\infty)$ and there are positive constants $m,M$ such that, given any $\phi\in S$, there exists $t_0=t_0(\phi)$  for which
$$m\le x_i(t,0,\phi)\le M,\q  1\le i\le n,\, t\ge t_0.$$
Hereafter,  unless otherwise stated, the notions of persistence, uniform persistence  and permanence  always refer to the choice of $S=C_0$ as the set of admissible initial conditions.
\end{defn}

Observe that a {\it linear homogeneous} DDE system is uniformly  persistent (in $C_0$) if and only if all components of all solutions with initial conditions in $C_0$ tend to $\infty$ as $t\to \infty$.  The next result concerns the uniform persistence of \eqref{7}.

\begin{prop}\label{prop1} Assume (h1), (h3), and  that  there exist  vectors $v\gg 0,\eta\gg 0$ such that
\begin{equation}\label{3.6}
M(t)v\ge \eta \q {\rm for\ large} \q t>0.
\end{equation}
 Then all solutions $y(t)$ of \eqref{7}
with initial conditions in $C_0$ satisfy $\lim_{t\to \infty} y_i(t)=\infty,\, i=1,\dots,n$.
\end{prop}

\begin{proof}
For $\phi\in C_0, t_0\ge \tau$, we have   $y_t=y_t(t_0,\phi)\in int\, C^+$ for $t\ge t_0$, thus $y_\tau\ge cv$ for some small $c>0$. System \eqref{7} is linear and cooperative, with $[v,\infty)$ forward invariant for $t$ on the  interval $[T_0,\infty)$ if $M(t)v\ge  0$ for $t\ge T_0$. To simplify the exposition, as before we take $v={\bf 1}$.
We only need to show that all components $u_i(t)$ of  the solution
$u(t):=y(t,T_0,{\bf 1})$ satisfy $\lim_{t\to\infty}u_i(t)=\infty$.

For $j\in\{1,\dots,n\}$, let $c_j=\liminf_{t\to\infty} x_j(t)\in[1,\infty]$. Suppose that $c_j<\infty$ for some $j\in\{1,\dots,n\}$, and take $c_i=\min_j c_j$, for the natural ordering in $(0,\infty]$.
Then, there is a sequence $t_k\to\infty$  such that $u_i(t_k)\to c_i, u_i'(t_k)\to 0$.  On the other hand, from \eqref{3.6}   there are $\eta_i, T_1>0$  such that
$\be_i(t)-d_i(t)+\sum_{j\ne i}a_{ij}(t)\ge \eta_i>0, t\, \ge T_1$.
 For any small $\vare>0$ and $k$ sufficiently large, we obtain
\begin{equation*}
\begin{split}
 u_i'(t_k)&\ge -d_i(t_k)u_i(t_k)+(c_i-\vare)\Big (\sum_{j\ne i} a_{ij}(t_k)+\be_{i}(t_k)\Big )\\
 &\ge d_i(t_k) [- u_i(t_k)+(c_i-\vare)]+(c_i-\vare)\eta_i,
\end{split}
\end{equation*}
and therefore $0\ge c_i\eta_i>0,$
which is not possible. This ends the proof.
\end{proof}

For dissipative systems \eqref{3} with nonlinearities satisfying (h4), the above criterion for the uniform persistence of the linearization at zero also provides  a criterion for its uniform persistence.  This is stated in the main theorem of this section, given below. For a relevant extension, see Theorem \ref{thm4++}.

\begin{thm}\label{thm4}  Assume (h1)-(h4), and suppose that  there exist $v\gg 0,\eta\gg 0$ such that \eqref{3.6}  is satisfied. Then \eqref{3}
is uniformly persistent, and thus permanent.
\end{thm}

\begin{proof}  After effecting a scaling of the variables, we take   $v=(1,\dots,1)={\bf 1}$ in condition  \eqref{3.6},  thus there exist constants $\eta_i>0\, ( i=1,\dots,n)$ such that, for some $T_0$,
$$\be_i(t)\ge d_i(t)-\sum_{j\ne i}a_{ij}(t)+\eta_i,\q t\ge T_0.$$
On the other hand, $d_i(t)-\sum_{j\ne i}a_{ij}(t)\le \ol d_i:=\sup _{t\ge T_0}d_i(t)$, and with $1<\al_i<1+\eta_i/\ol d_i$ we obtain
\begin{equation}\label{3.7}
\al_i^{-1}\be_i(t)-d_i(t)+\sum_{j\ne i}a_{ij}(t)>0,\q {\rm for}\q t\ge T_0,  i=1,\dots,n.
\end{equation}

 For $h_i^-$ as in (h4),
 we can choose
  $L>m>0$  such that the uniform estimate \eqref{3.5} holds, $h_i^-(m)=\min_{x\in [m,L]}h_i^-(x)$, with
  $(h_i^-)'(x)>0$ and $\al_i^{-1}x<h_i^-(x)$ for $x\in(0,m]$ and  all $i$.

Consider the auxiliary  {\it cooperative} system
\begin{equation}\label{05}
\begin{split}
x_i'(t)&=-d_i(t)x_i(t)+\sum_{j=1,j\ne i}^n a_{ij}(t)x_j(t)+\sum_{k=1}^m \be_{ik}(t)    H_i(x_i(t-\tau_{ik}(t)))\\
&=:F_i(t,x_t),\q  i=1,\dots,n,\ t\ge 0,
\end{split}
\end{equation}
where $H_i(x)=h_i^-(x)$ if $0\le x\le m$, $H_i(x)=h_i^-(m)$ if $x\ge m$.

For $x(t)$ a positive solution of \eqref{3},   for $t>0$ sufficiently large and $1\le i\le n$, we have $x_i(t)\le L$ and $h_{ik}(t,x_i(t-\tau_{ik}(t)))\ge H_i(x_i(t-\tau_{ik}(t)))$. Therefore, if \eqref{05} is uniformly persistent, then \eqref{3} is uniformly persistent as well.

Now,  we consider any solution $x(t)=x(t,T_0,\phi,F)$  of \eqref{05} with  $x_{T_0}=\phi\in C_0$ (where $T_0$ is as in \eqref{3.7}), and claim that
$$
\liminf_{t\to\infty} x_i(t)\ge m,\q 1\le i\le n.
$$
In fact, we shall show that there exists $T\ge 0$ such that
 \begin{equation}\label{3.9}
  x_i(t)\ge m\q {\rm for}\q  t\ge T,1\le i\le n.
\end{equation}
The proof, inspired by some arguments in \cite{FariaRost}, is divided into several steps.

\med

{\it Step 1}.  We  prove that if $\displaystyle \min\{ x_j(t):1\le j\le n,t\in [T,T+\tau]\}\ge m $ for some $T\ge T_0$, then $ x_j(t)\ge m $ for all $t\ge T$ and $j=1,\dots, n$.

\med

Assume that $x_j(t)\ge m$ for $t\in [T,T+\tau]$ and $j=1,\dots, n$. Let $t_0\in [T+\tau, T+2\tau]$ and $i\in \{1,\dots, n\}$ such that $x_i(t_0)=\min\{ x_j(t):1\le j\le n,t\in [T+\tau, T+2\tau]\}$.

If $x_i(t_0)<m$, we have 
$$ 0\ge x_i'(t_0)=-d_i(t_0)x_i(t_0)+\sum_{j=1,j\ne i}^n a_{ij} (t_0)x_j(t_0)+\sum_{k=1}^m \be_{ik} (t)H_i(x_i(t_0-\tau_{ik}(t_0))).
$$
Note that $x_i(t_0-\tau_{ik}(t_0))\ge m$ if $t_0-\tau_{ik}(t_0)\in [T,T+\tau]$ and $x_i(t_0-\tau_{ik}(t_0))\ge x_i(t_0)$ if $t_0-\tau_{ik}(t_0)\in [T+\tau, t_0]$, hence $H_i(x_i(t_0-\tau_{ik}(t_0))) \ge H_i(x_i(t_0))$. From \eqref{3.7} and the definition of $m$ we obtain
  \begin{equation}\label{06}
\begin{split}
0&\ge \left (-d_i(t_0)+\sum_{j=1}^n a_{ij}(t_0)\right ) x_i(t_0) +\be_i(t_0)H_i(x_i(t_0)) \\
&\ge \left (-d_i(t_0)+\sum_{j=1}^n a_{ij}(t_0)+\al_i^{-1}\be_i(t_0)\right )x_i(t_0)>0,
\end{split}
\end{equation}
which is not possible. Thus, $x_i(t_0)\ge m$. By iteration, this proves Step 1.

\med

{\it Step 2}.
Next,  for any $ T_0$ as in \eqref{3.7} and  $s_0:=\min\{ x_j(t):1\le j\le n,t\in [T_0,T_0+\tau]\}$, we shall show the estimate
$$\min_j\min_{t\in[T_0+\tau,T_0+2\tau]} x_j(t)\ge s_1,$$
where
$$s_1:=\min \Big\{ m,  \min_j \Big(\al_j H_j(s_0)\Big)\Big\}.$$

To simplify the exposition,  take $T_0=0$. In this way, we denote $s_0:=\min\{ x_j(t):1\le j\le n,t\in [0,\tau]\}>0$.
If $s_0\ge m$, from Step 1 the proof is complete.
Now, consider the case $s_0<m$.
By the definition of $m$, $h_j^-(s_0)  \al_j=H_j(s_0)  \al_j>s_0$ for all $j$, thus $s_1>s_0$. We claim that
 \begin{equation}\label{07}
 \min_j\min_{t\in[\tau,2\tau]} x_j(t)\ge s_1.
\end{equation}
Otherwise, there are $t_1\in [\tau, 2\tau]$ and $i\in\{ 1,\dots,n\}$ such that $x_i(t_1)<s_1$ and $x_j(t)\ge x_i(t_1)$ for all $t\in [\tau,t_1]$ and $j\in\{ 1,\dots,n\}$.


Since $x_i(t_1-\tau_{ik}(t_1))\ge \min \{s_0,x_i(t_1)\}$, we have
$H_i(x_i(t_1-\tau_{ik}(t_1)))\ge \min \{ H_i(s_0),H_i(x_i(t_1))\}$. We now consider two cases separately.

If $s_0\ge x_i(t_1)$, then $H_i(s_0)\ge H_i(x_i(t_1))$ and we get \eqref{06}  with $t_0$ replaced by $t_1$, thus a contradiction.

If $s_0< x_i(t_1)$, then $H_i(s_0)< H_i(x_i(t_1))$. Since $x_i(t_1)<s_1\le \al_iH_i(s_0)$, we derive
 \begin{equation*}
\begin{split}
0\ge x_i'(t_1)&\ge \left(-d_i(t_1)+\sum_{j=1}^n a_{ij}(t_1)\right) x_i(t_1)+\be_i (t_1)H_i(s_0) \\
&> \left(-d_i(t_1)+\sum_{j=1}^n a_{ij}(t_1)+\be_i (t_1)\al_i^{-1}\right)x_i(t_1)> 0,
\end{split}
\end{equation*}
which is again a contradiction. This proves the estimate \eqref{07}.

\med

{\it Step 3}.
Now, we define by recurrence the sequence
$$s_{k+1}=\min \Big\{ m,  \min_j \Big(\al_j H_j(s_k)\Big)\Big\},\q k\in \N_0.$$
If $s_k=m$ for some $k\in\N$, \eqref{3.9} follows by Steps 1 and 2.  In this case,
$\al_j H_j(s_k)> m,$
 hence $s_p=m$ for all $p>k$. If $s_k<m$ for all $k$, $(s_k)$ is strictly increasing, because
 $$s_{k+1}= \min_j \Big(\al_j H_j(s_k)\Big)>s_k.
 $$
 For $s^*=\lim s_k$,   from the definition of $m$ we derive
 $$0<s^*\le m \q {\rm and}\q s^*\ge \min_j \al_jH_j(s^*)> s^*,$$
which is not possible. The proof is complete.
\end{proof}

\begin{rmk}\label{rmk3.1} {\rm We observe that assumptions (h2) and \eqref{3.6} are satisfied if $\liminf_{t\to \infty}\be_i(t)>0$ and
\begin{equation}\label{09}
\ga_i\ge \frac{\be_i(t)v_i}{d_i(t)v_i-\sum_{j\ne i}a_{ij}(t)v_j}\ge \al_i>1,\q {\rm for}\q t\ge T_0,  \ i=1,\dots,n,
\end{equation}
for some vector $v=(v_1,\dots,v_n)\gg 0$ and constants $\al_i,\ga_i$.}\end{rmk}

\begin{exmp} {\rm Consider the system
 \begin{equation}\label{3.27}
x_i'(t)=-d_i(t)x_i(t)+\sum_{j=1,j\ne i}^n a_{ij}(t)x_j(t)+ \be_{i}(t)    \frac{x_i(t-\tau_{i}(t)))}{1+c_i(t)x_i^{\al_i}(t-\tau_{i}(t))},\  i=1,\dots,n,\ t\ge 0,
\end{equation}
where $\al_i\ge 1$, $d_i(t)>0, a_{ij}(t)\ge 0 ,c_i(t),\be_i(t),\tau_i(t)\ge 0$  are continuous and bounded and $0<c_i^-\le c_i(t)\le c_i^+, \be_i(t)\ge \be_i^->0$, $t\ge 0,i=1,\dots,n$. For  $h_i(t,x)=x(1+c_i(t) x^{\al_i})^{-1}, h_i^\pm(x)=x(1+c_i^\mp x^{\al_i})^{-1}$, we have $h_i^-(x)\le h_i(t,x)\le h_i^+(x)$, $h_i^\pm(0)=0,(h_i^\pm)'(0)=1$ and $0<h_i^-(x)\le h_i^+(x)<x$ for $x>0$. For each vector $v=(v_1,\dots,v_n)\gg 0$, define
$$l_i(v)=\liminf_{t\to\infty}{{\be_i(t)v_i}\over {d_i(t)v_i-\sum_{j\ne i}a_{ij}(t)v_j}},\ L_i(v)=\limsup_{t\to\infty}{{\be_i(t)v_i}\over {d_i(t)v_i-\sum_{j\ne i}a_{ij}(t)v_j}},\q i=1,\dots,n.$$
From Theorem \ref{thm3} and  Remark \ref{rmk3.1}, the zero solution of  \eqref{3.27} is GAS if there exists $v\gg 0$ such that $L_i(v)<1$ for all $i$,  whereas  \eqref{3.27} is permanent
if there exists $v\gg 0$ such that $1<l_i(v), L_i(v)<\infty$ for all $i$. }
\end{exmp}
\par
A careful reading of the proof above leads to several generalizations. First, it is clear that, in the statement of Theorem \ref{thm4},  hypothesis (h2)  can actually be replaced   by the dissipativeness of the system. Having this in mind, one also sees that the same arguments apply to dissipative systems more general than \eqref{3}, where, in each equation $i$, the instantaneous terms $a_{ij}(t) x_j(t)$ are replaced by linear delayed terms and the nonlinear terms are as in \eqref{0003}. This is expressed in the next theorem.


\begin{thm}\label{thm4++} Consider a non-autonomous  system
 of one of the forms
\begin{equation}\label{3++}
\begin{split}
x_i'(t)&=-d_i(t)x_i(t)+\sum_{j=1,j\ne i}^n  a_{ij}(t)\int_{-\tau}^0x_j(t+s)\, d_s\nu_{ij}(t,s)\\
&+\sum_{k=1}^m \be_{ik}(t)  \int_{-\tau}^0  h_{ik}(s,x_i(t+s))\,  d_s\eta_{ik}(t,s),\  i=1,\dots,n,\ t\ge 0,
\end{split}
\end{equation}
\begin{equation}\label{3+++}
\begin{split}
x_i'(t)&=-d_i(t)x_i(t)+\sum_{j=1,j\ne i}^n  a_{ij}(t)\int_{-\tau}^0x_j(t+s)\, d_s\nu_{ij}(t,s)\\
&+\sum_{k=1}^m \be_{ik}(t)   h_{ik}\left(t, \int_{-\tau}^0x_i(t+s)\,  d_s\eta_{ik}(t,s)\right),\  i=1,\dots,n,\ t\ge 0,
\end{split}
\end{equation}
where: $d_i(t), a_{ij}(t), \be_{ik}(t), h_{ik}(t,x)$ satisfy  (h1), (h3)  and (h4); the measurable functions $\nu_{ij},\eta_{ik}:[0,\infty)\times[-\tau,0]\to \R$ are continuous from the left in $s$, $\nu_{ij}(t,\cdot), \eta_{ik}(t,\cdot)$ are non-decreasing and normalized so that
$$\int_{-\tau} ^0 d_s\nu_{ij}(t,s)= \int_{-\tau} ^0 d_s\eta_{ik}(t,s)= 1 ,\q i,j=1,\dots,n,k=1,\dots,m,t\ge 0.$$
Assume also that \eqref{3++}, or \eqref{3+++}, is dissipative and that   \eqref{3.6}  is satisfied for some $v,\eta\gg 0$. Then, the system
is uniformly persistent.
\end{thm}

 We now apply the previous theorems to the Nicholson system \eqref{04}.

\begin{thm}\label{thm5} Consider system \eqref{04},   where  $d_i,a_{ij}, \be_{ik}, c_{ik}, \tau_{ik}:[0,\infty)\to [0,\infty)$ are  continuous and bounded,
 with $d_i(t), \be_i(t)=\sum_{k=1}^m \be_{ik}(t)$ strictly positive and $c_{ik}(t)\ge c_i>0$ on $[0,\infty)$, for all $i,j,k$.  With the notation in  \eqref{3.2}, assume that there exist a vector $u\gg 0$  and $\de>0, T_0\ge 0$ such that $[D(t)-A(t)]u\ge \de{\bf 1}$ for $t\ge T_0$. Then:
 \vskip 0cm
 (i)  Eq. \eqref{04} is dissipative. \vskip 0cm
 (ii) If there exist a vector $v\gg 0$   and $T_0\ge 0$ such that $M(t)v\le 0$ for $t\ge T_0$ and either $\liminf_{t\to\infty} \be_i(t)\gg 0$ or  $\limsup_{t\to\infty} (M(t)v)_i\ll 0\, (1\le i\le n)$,  the zero solution of \eqref{04} is GAS.
  \vskip 0cm
 (iii) If there exist  vectors $v\gg 0,\eta\gg 0$ such that $M(t)v\ge  \eta$ for $t\ge T_0$, then \eqref{04} is  permanent.
\end{thm}

%


\begin{exmp}  {\rm Consider the planar system
 \begin{equation}\label{3.28}
 \begin{split}
x_1'(t)&=-(1+\cos ^2t)x_1(t)+\ga_1(1+\sin^2t)x_2(t)+
\sum_{j=1}^{m} (\be_{1j}+f_{1j}(t))x_1(t-\tau_{1j}(t))e^{-c_{1j}(t)x_1(t-\tau_{1j}(t))},\\
x_2'(t)&=-(1+\sin ^2t)x_2(t)+\ga_2(1+\cos^2t)x_1(t)+
\sum_{j=1}^{m} (\be_{2j}+f_{2j}(t))x_2(t-\tau_{2j}(t))e^{-c_{2j}(t)x_2(t-\tau_{2j}(t))},
\end{split}
\end{equation}
where $\ga_i>0,\be_{ij}\ge 0$ with $\be_i:=\sum_{j=1}^m\be_{ij}>0$, all the functions $f_{ij}(t), c_{ij}(t), \tau_{ij}(t)$ are continuous, nonnegative and bounded on $[0,\infty)$, with $c_{ij}(t)$ bounded below by positive constants, for $t\ge 0,i=1,2, j=1,\dots,m$. Write $f_i(t):=\sum_{j=1}^mf_{ij}(t)$, let $ f_i^-,f_i^+$ be such that $0\le f_i^-\le f_i(t)\le f_i^+$ for $t\ge 0$ and denote $\be_i^-=\be_i+f_i^-,\be_i^+=\be_i+f_i^+, \, i=1,2$. With the notation in \eqref{3.2}, we have
 \begin{equation*}
 \begin{split}
D(t)-A(t)&=\left [\begin{matrix}
1+\cos^2t&-\ga_1(1+\sin^2t)\\
-\ga_2(1+\cos^2t)&1+\sin^2t
\end{matrix}\right ],\\
M(t)&=\left [\begin{matrix}
\be_1+f_1(t)-(1+\cos^2t)&\ga_1(1+\sin^2t)\\
\ga_2(1+\cos^2t)&\be_2+f_2(t)-(1+\sin^2t)\\
\end{matrix}\right ].
\end{split}
\end{equation*}

Consider a vector $u=(1,u_2)$ with $u_2>0$, and write $[D(t)-A(t)]u=\left [\begin{matrix} \de_1(t)\\ \de_2(t)\end{matrix}\right ]$. Since $\min\de_1(t)=1-2u_2\ga_1, \min\de_2(t)=u_2-2\ga_2$, if
 \begin{equation}\label{3.30a}
 4\ga_1\ga_2\le 1
 \end{equation}
 we can find $u_2$ such that $2\ga_2\le u_2\le (2\ga_1)^{-1}$, implying that $[D(t)-A(t)]u\ge 0$ for $t\in\R$. On the other hand, $\de_1(\frac{\pi}{4})=\frac{3}{2}(1-u_2\ga_1)>0, \de_2(\frac{\pi}{4})=\frac{3}{2}(-\ga_2+u_2)>0$. By Theorem \ref{thm20}, we conclude that the ODE $x'=-[D(t)-A(t)]x$ is globally  exponentially stable.

  We now look for a vector $v=(1,v_2)\gg 0$ such that  $M(t)v\ge \eta\gg 0$. Write $M(t)v=\left [\begin{matrix} m_1(t)\\ m_2(t)\end{matrix}\right ]$ and observe that
 $m_1(t)\ge \eta_1:=\be_1^--2+v_2\ga_1,\, m_2(t)\ge\eta_2:= v_2(\be_2^--2)+\ga_2,\, t\ge 0.$
 Now assume that:
 \begin{equation}\label{3.30b}
{\rm either}\q \be_1^-\ge 2\q {\rm or}\q \be_2^-\ge 2\q {\rm or}\q
 (2-\be_1^-)(2-\be_2^-)<\ga_1\ga_2. \end{equation}
One easily verifies that: (i) if either $\be_1^-\ge 2$ or $\be_2^-\ge 2$, one can find $v_2>0$ such that
$\eta_1>0,\eta_2>0$; (ii)  if $\be_i^-<2$ for $i=1,2$, and $ (2-\be_1^-)(2-\be_2^-)<\ga_1\ga_2$,
for any $v_2$ such that  $(2-\be_1)\ga_1^{-1}<v_2<\ga_2(2-\be_2)^{-1}$ we have $M(t)v\ge \eta=(\eta_1,\eta_2)\gg 0$. From Theorem \ref{thm5}.(iii), conditions \eqref{3.30a}-\eqref{3.30b} imply that  \eqref{3.28} is permanent.

As an illustration, with  $m=1$ and  $\be_1=2,\ga_1=1=4\ga_2,f_1(t)=f_2(t)=0$, we conclude that \begin{equation}\label{3.29}
 \begin{split}
x_1'(t)&=-(1+\cos ^2t)x_1(t)+(1+\sin^2t)x_2(t)+
2x_1(t-\tau_{1j}(t))e^{-c_{1j}(t)x_1(t-\tau_{1j}(t))}\\
x_2'(t)&=-(1+\sin ^2t)x_2(t)+\frac{1}{4}(1+\cos^2t)x_1(t)+
\be_2 x_2(t-\tau_{2j}(t))e^{-c_{2j}(t)x_2(t-\tau_{2j}(t))}\
\end{split}
\end{equation}
is permanent for any $\be_2>0$.

On reverse, if  $\be_i^+<1, \, i=1,2$, for a positive vector $v=(1,v_2)$  we obtain $M(t)v=\left [\begin{matrix} m_1(t)\\ m_2(t)\end{matrix}\right ]$ with
 $m_1(t)\le \eta_1:=\be_1^+-1+2v_2\ga_1,\ m_2(t)\le\eta_2:= v_2(\be_2^+-1)+2\ga_2,\ t\ge 0.$
 At this point, assume
 \begin{equation}\label {3.30c}
 \be_i^+<1,\ i=1,2 \q {\rm and}\q 4\ga_1\ga_2\le (1-\be_1^+)(1-\be_2^+).
 \end{equation}
 Thus, choosing $v_2$ such that $2\ga_2(1-\be_2^+)^{-1}\le v_2\le (2\ga_1)^{-1}(1-\be_1^+)$ we obtain $M(t)v\le 0$ for all $t\ge 0$. From Theorem \ref{thm5}.(ii), conditions \eqref{3.30c} imply that the trivial solution of \eqref{3.28} is GAS. In particular, this is the case of the zero solution of
 \begin{equation}\label{3.29}
 \begin{split}
x_1'(t)&=-(1+\cos ^2t)x_1(t)+\frac{1}{4}(1+\sin^2t)x_2(t)+
\frac{1}{2}x_1(t-\tau_{1j}(t))e^{-c_{1j}(t)x_1(t-\tau_{1j}(t))}\\
x_2'(t)&=-(1+\sin ^2t)x_2(t)+\frac{1}{4}(1+\cos^2t)x_1(t)+
\be_2 x_2(t-\tau_{2j}(t))e^{-c_{2j}(t)x_2(t-\tau_{2j}(t))}\
\end{split}
\end{equation}
for any $0<\be_2\le\frac{1}{2}$.}
\end{exmp}

\begin{rmk}\label{rmk3.3} {\rm In recent years, some attention has been given to  Nicholson's blowflies equations and systems with {\it  harvesting}. For the $n$-dimensional case, such systems are obtained by adding linear harvesting terms with delays to \eqref{04}, so that it becomes:
\begin{equation}\label{3.31}
\begin{split}
x_i'(t)=&-d_i(t)x_i(t)+\sum_{j=1,j\ne i}^n a_{ij}(t)x_j(t)\\
&+\sum_{k=1}^m \be_{ik}(t)    x_i(t-\tau_{ik}(t))e^{-c_{ik}(t)x_i(t-\tau_{ik}(t))}
-H_i(t)x_i(t-\sigma_i(t)),
\q  i=1,\dots,n,
\end{split}
\end{equation}
where the new coefficients $H_i(t)$ and delays $\sigma_i(t)$ are continuous, nonnegative and bounded. For the scalar case of \eqref{3.31},
 Liu \cite{Liu13} studied both  the global exponential stability of the zero solution and  the permanence. The almost periodic scalar case of \eqref{3.31} was studied in \cite{Zhang}, and the $n$-dimensional  case in
  \cite{Wang13},  where the authors  established criteria for
 the existence and global exponential stability of a positive almost periodic solution by using properties of almost periodic functions and Lyapunov functionals. See also \cite{Zhou} for a periodic  system \eqref{3.31} with $n=2$. From the proof of Theorem \ref{thm3}, we deduce that Theorem \ref{thm5}.(ii), on the global asymptotic stability of the zero solution, applies to \eqref{04} replaced by \eqref{3.31},   without any changes. However, the result on permanence  in Theorem \ref{thm5}.(iii) does not carry over to \eqref{3.31}. An interesting open problem is to generalize our results, and find sufficient conditions for the permanence of \eqref{3.31}.}

\end{rmk}

\begin{rmk}\label{rmk3.3}{\rm In \cite{Faria15}, Faria studied the persistence and permanence of a class of {\it cooperative}  DDEs with possible {\it infinite} delay of the form $x_i'(t)=F_i(x_t)-x_i(t)G_i(x_t),\, 1\le i\le n$. By using properties of cooperative  systems, it was shown that, under some additional conditions, all positive solutions  are bounded below and above  by  positive equilibria, which in particular proves the permanence.  The persistence and permanence for the  non-autonomous system  $x_i'(t)=F_i(t,x_t)-x_i(t)G_i(t,x_t),\, 1\le i\le n$, was also addressed in  \cite{Faria15}
   by comparing it above and below with autonomous cooperative systems. Although
   the  basic idea  is    similar (comparison of solutions with solutions of cooperative systems), the results and  techniques in  \cite{Faria15} do not apply to the study of systems \eqref{3}: not only does \eqref{3}  not have the above form, but  the nonlinearities  $h_{ik}(t,x)$ are in general non-monotone on the second variable. On the other hand, this remark raises another interesting open problem: how to extend  the  results about permanence   in this paper  to systems with {\it infinite delay}, since it is clear that the proof of Theorem \ref{thm4} does not work for the infinite delay case.}
    \end{rmk}

\section{Sharp criteria for systems with autonomous coefficients}
\setcounter{equation}{0}

The   case of an autonomous system \eqref{3}, or of \eqref{3} with constant coefficients but time-varying delays,  is particularly important in applications.
For these situations, the matrices $A,B,D,M$ in \eqref{3.2} are autonomous, and their properties play an important role in the analysis of the asymptotic behaviour  of solutions.
For the sake of completeness and convenience of the reader,  some elements from matrix theory will be recalled here. We start with some definitions.
\begin{defn}\label{def4.1} Let $N=[n_{ij}]$ be a square matrix. The matrix
  $N$ is said to be   {\bf reducible} if there is a simultaneous permutation of rows and columns that brings $N$ to the form
  $$ \left[\begin{array}{cc}
   N_{11}&0\\ N_{21}&N_{22}\end{array}\right],$$
  with $N_{11}$ and $N_{22}$   square matrices;
     $N$ is an {\bf irreducible matrix} if it is not reducible.
For $N$ with nonpositive off-diagonal entries   (i.e., $n_{ij}\le 0$ for $i\ne j$), $N$ is said to be a {\bf non-singular M-matrix} 
if   all its eigenvalues have positive 
real parts. We say that $N$ is a
{\bf cooperative  matrix} if it has  nonnegative off-diagonal entries (i.e., $n_{ij}\ge 0$ for $i\ne j$).
\end{defn}

The reader should be aware  that many authors use the term {\it M-matrix} with the above meaning of the term  {\it non-singular M-matrix}.
For alternative definitions and  properties of M-matrices, see \cite{Fiedler}.
Namely,  it is important to remark that,  for a square matrix $N$ with nonpositive off-diagonal entries, $N$  is a
non-singular M-matrix if and only if
there exists a vector $u\gg 0$ such that $Nu\gg 0$.

System  \eqref{3} with constant coefficients becomes
\begin{equation}\label{3.14}
x_i'(t)=-d_ix_i(t)+\sum_{j=1,j\ne i}^n a_{ij}x_j(t)+\sum_{k=1}^m \be_{ik}    h_{ik}(x_i(t-\tau_{ik}(t))),\  i=1,\dots,n,\ t\ge 0,
\end{equation}
 and hypothesis (h4) translates simply as
\begin{itemize}
\item[(h4*)] $h_{ik}: [0,\infty)\to [0,\infty)$ are  bounded,  locally Lipschitzian and continuously differentiable on a vicinity of $0^+$,  with $h_{ik}(0)=0, h_{ik}'(0)=1$ and  $h_{ik}(x)>0$ for $x>0$,   $i\in \{1,\dots,n\}, k\in\{1,\dots,m\}$.
\end{itemize}

For \eqref{3.14}, the results in the previous section are summed up in the following theorem:


 \begin{thm}\label{thm4.1}  Consider system \eqref{3.14},
where $d_i>0$, $a_{ij}\ge 0, \be_{ik}\ge 0$ with $\be_i:=\sum_{k=1}^m \be_{ik}>0$,
$h_{ik}, \tau_{ik}:[0,\infty)\to [0,\infty)$  are continuous, with $ \tau_{ik}(t)$ uniformly bounded from above by some $\tau>0$,  $i,j=1,\dots,n, k=1,\dots,m$, and suppose that   (h4*) is satisfied.
Define the $n\times n$ matrices
\begin{equation}\label{3.15}
A=[a_{ij}],\ B=diag\, (\be_1,\dots,\be_n),\ D=diag\, (d_1,\dots,d_n),\ M=B-D+A,
\end{equation}
where $a_{ii}:=0\ (1\le i\le n)$, and assume that $D-A$ is a non-singular M-matrix.
  Then:
 \vskip 0cm
 (i) \eqref{3.14} is dissipative;\vskip 0cm
 (ii) If   in addition $h_i^+(x):=\max_{1\le k\le m} h_{ik}(x)<x$ for $x>0,   i=1,\dots,n,$ and there exists a vector $v\gg 0$    such that $Mv\le 0$,  the trivial solution of \eqref{3.14} is GAS;
  \vskip 0cm
 (iii) If there exists a  vector $v\gg 0$ such that $Mv\gg 0$,  \eqref{3.14} is  permanent.
\end{thm}

For  an $n\times n$ matrix $N$,  the {\it spectral bound} or {\it stability modulus} $s(N)$ is defined by
$$s(N)=\{ Re\, \la:\la\in \sigma(N)\},$$
where $\sigma(N)$ denotes the spectrum of $N$. For a
cooperative and irreducible matrix $N$,  it is well-known that the spectral bound $s(N)$ is a (simple)  eigenvalue, with  a  strictly positive associated eigenvector, see Appendix A.5 of \cite{SmithThi}; moreover, 
$s(N)>0$ if and only if there exists a strictly positive vector $v\in\R^n$ with $Nv\gg 0$
 \cite{FariaRost}.  Thus,  a threshold  criterion of permanence versus extinction is obtained from Theorem \ref{thm4.1} when $A$ is  an irreducible matrix.

 \begin{cor}\label{cor3.2} Assume  all the general hypotheses of Theorem \ref{thm4.1} (including (h4*) and that $D-A$ is a non-singular M-matrix) are satisfied. Further assume that $A$ is irreducible and $h_i^+(x)<x$ for $x>0,   i=1,\dots,n$. Then:
 (i) if $s(M)\le 0$,  the trivial solution of \eqref{3.14} is GAS;
  (ii) if $s(M)>0$,  \eqref{3.14} is  permanent.
\end{cor}

This threshold  criterion is not valid, in general, when $A$ (and therefore $M$ as well) is  reducible. Our next task is to replace the assumptions in
Theorem \ref{thm4.1} by sharp conditions for extinction versus permanence when $M$ is reducible. We emphasize that usually the case of a reducible community matrix $M$ is  not treated in the literature.

By an adequate simultaneous permutation of rows and columns, which amounts to a permutation of the variables in the original system \eqref{3.14},
we may suppose that the  $n\times n$-matrix $A$
has been transformed into the  triangular form
\begin{equation}\label{triangular}
A = \left[\begin{array}{cccc}
A_{11} & 0 & \ldots & 0 \\
A_{21} & A_{22} & \ldots& 0 \\
\vdots & \vdots  &\ddots & \vdots \\
A_{k1} & A_{k2} & \ldots& A_{kk}
\end{array}\right]\,,
\end{equation}
where the diagonal blocks $A_{11},\ldots, A_{kk}$ are square matrices  of size
$n_1,\ldots,n_k$ respectively, $n_1+\cdots+n_k=n$, and  are irreducible.
 Clearly,   $k=1$ if $A$ is
 irreducible. Observe that  a square  $n\times n$-matrix $A=[a_{ij}]$ is {\em
irreducible\/} if and only if  for any nonempty proper subset
$I\subset\{1,\ldots,n\}$ there are $i\in I$ and $j\in
\{1,\ldots,n\}\setminus I$ such that $a_{ij}\not=0$.

The next result extends Corollary \ref{cor3.2} and gives necessary and sufficient conditions for both the uniform persistence  and the global asymptotic stability of the zero solution of \eqref{3.14}, in the case of a reducible matrix $A$. The result for uniform persistence was inspired by~\cite{noos7}.

\begin{thm}\label{thmOS}
Consider system \eqref{3.14}
where $d_i>0$, $a_{ij}\geq 0$, $\beta_{ik}\geq 0$ with $\beta_i:=\sum_{k=1}^m\beta_{ik}>0$, $\tau_{ik}:[0,\infty)\to [0,\infty)$ are continuous and  bounded,   the functions $h_{ik}$ satisfy (h4*) with $h_i^+(x):=\max_{1\le k\le m} h_{ik}(x)<x$ for any $x> 0$, for $i,j=1,\ldots,n,\,k=1,\ldots,m$. Let $A,B,D,M$ be the matrices  defined in \eqref{3.15}.
 Assume that $D-A$  is a non-singular $M$-matrix.
Without loss of generality,  further assume that $A$ has the block lower triangular structure as in~\eqref{triangular}, with irreducible diagonal blocks $A_{11},\ldots,A_{kk}$, and denote by $M_{jj}$ the associated blocks in the matrix $M$, that is, $M_{jj}=B_{j}-D_{j}+A_{jj}$, with $B_j=diag(\beta_i)_{i\in I_j}$ and $D_j=diag(d_i)_{i\in I_j}$, where $I_j$ is  the set formed by
the $n_j$ indexes corresponding to the rows of the block $A_{jj}$, for each $j=1,\dots,k$; and $M_{ij}=A_{ij}$ for $1\leq j<i\leq k$.
Then:\vskip 0cm
(i) System \eqref{3.14} is uniformly persistent if and only if $s(M_{jj})>0$ for every index $j\in \{1,\ldots,k\}$ such that, except for the diagonal block $M_{jj}$, all the other blocks on the row are null.\vskip 0cm
(ii) The null solution of system \eqref{3.14} is GAS if and only if $s(M)\leq 0$.
\end{thm}

\begin{proof}  In the case of $A$ irreducible,  the results are given in  Corollary \ref{cor3.2}.  From now on, $A$ is assumed to be a reducible matrix with the triangular form~\eqref{triangular} with $k>1$.  We make a few remarks beforehand.

First, observe that the property of $D-A$ being a non-singular M-matrix is preserved under a simultaneous permutation of rows and columns (so that $D-A$ becomes $P(D-A)P^T$ for some orthogonal matrix $P$), therefore system  \eqref{3.14} is dissipative, and thus the uniform persistence in (i) can actually be  replaced by the permanence. Secondly, for each $j=1,\ldots,k$, we consider the lower dimensional system associated with the irreducible block $M_{jj}$,  formed by the $n_j$ equations
\begin{equation}\label{sistema j}
x_i'(t)=-d_i\,x_i(t) +\dps \sum_{p\in I_j,p\not= i}  a_{ip}\,x_p(t) + \sum_{q=1}^m \beta_{iq}\,h_{iq}(x_i(t-\tau_{iq}(t)))\,,\; i\in I_j\,,\; t\geq 0\,,
\end{equation}
and observe that it  satisfies all the hypotheses in Corollary \ref{cor3.2}, as $D_j-A_{jj}$ is a non-singular M-matrix as well.
Finally,  for any vector $v$ in $\R^n$ or any map $v$ taking values in $\R^n$, we introduce the notation $v^j=(v_i)_{i\in I_j}$ for each $j=1,\ldots,k$, so that $v=(v^1,\dots,v^k)$.
\med

(i) To simplify the writing, we may assume without loss of generality that the diagonal blocks in~\eqref{triangular} with all null blocks to their left (if any) are placed in the first rows. In other words, we assume that $\{A_{11},\ldots,A_{ll}\}$ are exactly the diagonal blocks with all other blocks on their row null, for some $1\leq l\leq k$. In this way, for $1\leq j\leq l$, system~\eqref{sistema j} is just a lower dimensional decoupled subsystem of system~\eqref{3.14}.
\par
Now, suppose that system \eqref{3.14} is uniformly persistent.  Then, for each  $j=1,\ldots,l$, system~\eqref{sistema j} naturally inherits the property of uniform persistence from the total system, and Corollary \ref{cor3.2}  implies that  $s(M_{jj})>0$ for any $j=1,\ldots,l$, so  this is a necessary condition.

Conversely, assume that $s(M_{jj})>0$ for any $j=1,\ldots,l$. Applying once more Corollary \ref{cor3.2}, we deduce that systems~\eqref{sistema j} are uniformly persistent for any $j=1,\ldots,l$. Therefore, there exists $m_0>0$ such that for any $\phi\in C_0$, $\liminf_{t\to\infty} x_i(t,0,\phi)\geq m_0$ for all $i\in I_1\cup \ldots \cup I_l$. At this point, if $l=k$ the proof is complete, whereas if $l<k$ we have to deal with the remaining components of the solution.
\par
We now consider the case $l<k$ and look at the components $x_i(t,0,\phi)$ for $i\in I_{l+1}$. The method here is twofold: first, since there is at least a non-null block to the left of $M_{l+1,l+1}$,  we will show that  one component   $x_{i_1}(t,0,\phi)$ ($i_1\in I_{l+1}$) of the solution eventually stays bounded away from $0$. Secondly, once we have raised one component in $I_{l+1}$, we recursively raise the rest of them, one by one, by applying the irreducible character of $M_{l+1,l+1}$.
\par
More precisely, as there is at least one non-null block to the left of $M_{l+1,l+1}$, there are indexes $i_1\in I_{l+1}$ and $j_1\in I_j$ for some $1\leq j\leq l$ such that $a_{i_1j_1}>0$. Now, for  an initial condition $\phi\in C_0$, there exists a $t_0=t_0(\phi)$ such that $x_i(t,0,\phi)\geq m_0$ for all $t\geq t_0$ and for all $i\in I_1\cup \ldots \cup I_l$. Therefore, for  $t\geq t_0$, $x_{i_1}'(t,0,\phi)\geq -d_{i_1}\,x_{i_1}(t,0,\phi) + a_{i_1j_1}\,m_0$. Now, we consider the scalar cooperative ODE
\begin{equation*}
y'(t)=-d_{i_1}\,y(t) + a_{i_1j_1}\,m_0\,,\; t\geq 0\,,
\end{equation*}
whose solution, for the previous time $t_0\geq 0$, is written as
$$
y(t,t_0,y(t_0))=y(t_0)\,e^{-d_{i_1}(t-t_0)}+\frac{a_{i_1j_1}\,m_0}{d_{i_1}}(1-e^{-d_{i_1}(t-t_0)})\,,
$$
so that  there exist $m_{1}>0$ and $t_{1}\geq t_0$ such that $y(t,t_0,y(t_0))\geq m_{1}$ for any $t\geq t_{1}$, provided that $y(t_0)\geq 0$.
The application of a standard argument of comparison of solutions permits to conclude  that
$x_{i_1}(t,0,\phi)\geq m_{1}$ for any $t\geq t_{1}$.
\par
If $I_{l+1}=\{i_1\}$, we are done with this block. If not, as $A_{l+1,l+1}$ is irreducible, there exists an index $i_2\in I_{l+1}\setminus\{i_1\}$ such that $a_{i_2i_1}>0$.
As before, we consider the scalar ODE
\begin{equation*}
y'(t)=-d_{i_2}\,y(t) + a_{i_2i_1}\,m_1\,,\; t\geq 0\,,
\end{equation*}
for which we find a constant $m_2>0$ and a  time $t_2\geq t_1$ such that if $t\geq t_2$, $y(t,0,y(t_1))\geq m_2$ for any $t\geq t_2$, independently of the value $y(t_1)\geq 0$. In a similar way, we conclude that the $i_2$th component of the solution of~\eqref{3.14} satisfies $x_{i_2}'(t,0,\phi)\geq -d_{i_2}\,x_{i_2}(t,0,\phi) + a_{i_2i_1}\,m_1$ for any $t\geq t_1$, and once more, by comparing solutions, we have $x_{i_2}(t,0,\phi)\geq m_2$ for any $t\geq t_2$.
\par
At this point, if $I_{l+1}=\{i_1,i_2\}$  we are finished with this block;  if not, as $A_{l+1,l+1}$ is irreducible, considering  $\{i_1,i_2\}$ and its complement  $I_{l+1}\setminus\{i_1,i_2\}$, we may affirm that there exist indexes $i_3\in I_{l+1}\setminus\{i_1,i_2\}$ and $j\in \{i_1,i_2\}$ such that $a_{i_3j}>0$; now, the argument to lift the component $x_{i_3}(t,0,\phi)$ is just the same as the one for $x_{i_2}(t,0,\phi)$.
\par
Iterating this procedure inside the irreducible block  $A_{l+1,l+1}$, we  conclude that there is a constant $m_0'=\min\{m_0,m_1,\ldots,m_{n_{l+1}}\}>0$ such that for any $\phi\in C_0$, there exists a $t_0'=t_0'(\phi)$ such that $x_i(t,0,\phi)\geq m_0'$ for all $t\geq t_0'$ and for all $i\in I_1\cup \ldots \cup I_{l}\cup I_{l+1}$.
\par
To finish,  note that the procedure for the remaining components of the solution, if any, is identical to the one just developed for the set of indexes $I_{l+1}$.

\med

(ii) Note that $s(M)=\max\{s(M_{11}),\ldots,s(M_{kk})\}$, so that $s(M)\leq 0$ if and only if  $s(M_{jj})\leq 0$ for $j=1,\ldots,k$.
Because of the triangular structure of $A$ in \eqref{triangular}, and
with  the previous notation for $\phi=(\phi^1,\dots,\phi^k)$, it is apparent that,  for $j=2,\ldots,k$, the ``faces"
$$
F_j= \{\phi=(\phi^1,\ldots,\phi^k)\in C^+\mid \phi^1=\ldots=\phi^{j-1}=0\}
$$
of the nonnegative cone $C^+$ are
positively invariant. In this way, for an initial condition $\phi\in F_{j}\ (2\le j\le k)$, the
solution remains in $F_{j}$, thus the component  $x^{j}(t,0,\phi)$  is a solution of the system  \eqref{sistema j}.

\smal

We first assume that the null solution of \eqref{3.14} is GAS in the nonnegative cone $C^+$.
For any $j=1,\ldots,k$ fixed,  we now show  that any solution $y(t,0,\phi^j)$ of system (4.4) with initial condition $\phi^j\in C^+([-\tau,0];\R^{n_j})$ has $\lim_{t\to\infty} y(t,0,\phi^j)=0$. This is clear for $j=1$, as system \eqref{sistema j} is a decoupled subsystem of system \eqref{3.14}. For $j>1$, just consider $\widetilde\phi\in C^+([-\tau,0];\R^{n})$ with $\widetilde\phi^j=\phi^j$ and $\widetilde\phi^i=0$ for any $i<j$, so that $\widetilde\phi\in F_j$. Then,  $x^j(t,0,\widetilde\phi)$ is a solution of system \eqref{sistema j}, thus $y(t,0,\phi^j)=x^j(t,0,\widetilde\phi)\to 0$ as $t\to\infty$, as we wanted.  With this behaviour for each $j$, the persistent case in Corollary \ref{cor3.2} is precluded, and then it must be $s(M_{jj})\leq 0$.



\smal

Conversely, assume that $s(M_{jj})\leq 0$, so that there exist vectors $v^j\in\R^{n_j}$, $v_j\gg 0$ such that $M_{jj}\,v^j\leq 0$ for  $j=1,\dots,k$. Without loss of generality, we suppose that $k=2$, so that
$M$ has the form
\begin{equation}\label{M2}
M=B- D+ \left[\begin{array}{cc}
 A_{11} & 0  \\
 A_{21} & A_{22}
\end{array}\right],
\end{equation}
where $A_{ii}$  are $n_i\times n_i$ irreducible blocks ($i=1,2$).
 The general case of $k$ blocks follows by iterating the procedure below.

With $M$ given by \eqref{M2},  \eqref{3.14} is equivalent to
\begin{equation}\label{systemM2}
\begin{split}
x_i'(t)&=-d_i\,x_i(t) +\dps \sum_{p\in I_1,p\not= i}  a_{ip}\,x_p(t) + \sum_{k=1}^m \beta_{ik}\,h_{ik}(x_i(t-\tau_{ik}(t))),\ i\in I_1,\, t\geq 0\\
x_i'(t)&=-d_i\,x_i(t) +\dps \sum_{p\in I_2,p\not= i}  a_{ip}\,x_p(t)+ \sum_{k=1}^m \beta_{ik}\,h_{ik}(x_i(t-\tau_{ik}(t))) +\sum_{p\in I_1} a_{ip} x_p(t),\  i\in I_2,\, t\geq 0\,.
\end{split}
\end{equation}

We first claim that the trivial solution of \eqref{3.14} is globally attractive.  Let $\phi\in C^+$, and write $x(t,0,\phi)=(x^1(t,0,\phi^1),x^2(t,0,\phi))$. (Recall that $x^1(t,0,\phi)$ is just the solution of system~\eqref{sistema j} for $j=1$ with initial condition $\phi^1$.) Corollary \ref{cor3.2} implies that $\lim_{t\to \infty}x^1(t,0,\phi^1)=0$. In particular, in
\eqref{systemM2} we have $q_i(t):= \sum_{p\in I_1} a_{ip} x_p(t,0,\phi^1)\to 0$ as $t\to\infty$ for any $i\in I_2$. At this point,
the proof of $\lim_{t\to \infty}x^2(t,0,\phi)=0$ is obtained by simply  repeating the
 argument used in the proof of Theorem \ref{thm3}  applied to the  second system in \eqref{systemM2}. Details are omitted.

\smal

It remains to prove the stability of  the null solution of \eqref{systemM2}.    For a given $\phi=(\phi^1,\phi^2)\in C^+$, as before we denote the solution of \eqref{3.14} by $x(t,0,\phi)=(x^1(t,0,\phi^1),x^2(t,0,\phi))$.

From the assumptions on $h_{ik}$, for $i\in I_2$, we construct  maps $\widetilde h_i:[0,\infty)\to[0,\infty)$
satisfying the following conditions: $\widetilde h_i$ are continuous, bounded, nondecreasing, equal to $h_i^+$ on a right neighbourhood of $0$, and such that $h_i^+(x)\leq \widetilde h_i(x)<x$  for all $x> 0$.
Now, we consider an $n$-dimensional system, whose first $n_1$ equations are given by \eqref{sistema j} with $j=1$ (as in \eqref{systemM2}) and the last equations given by the $n_2$-dimensional system
\begin{equation}\label{4.8}
y_i'(t)=-d_i\,y_i(t) +\dps \sum_{p\in I_2,p\not= i}  a_{ip}\,y_p(t)+ \sum_{k=1}^m \beta_{ik}\,\widetilde h_i(y_i(t-\tau_{ik}(t)))+\sum_{p\in I_1}  a_{ip}\,x_p(t),\ i\in I_2, t\ge 0,
\end{equation}
 written for short as $y_i'(t)=\widetilde f_i(t,y_t)$. Since  \eqref{4.8} is cooperative,  a comparison of solutions leads to $x^2(t,0,\phi)\leq y(t,0,\phi, \widetilde f)$ for  $t\geq 0$.

\smal

Fix any $\varepsilon>0$. Let $\widetilde\varepsilon>0$ be sufficiently small so that $\beta_i(\varepsilon v_i^2- \widetilde h_i(\varepsilon v_i^2))\ge \widetilde\varepsilon\, v_i^2,\ i\in I_2$.
Of course,  Corollary \ref{cor3.2} (or Theorem \ref{thm3}) yields  the stability of the null solution for the first  system in \eqref{systemM2},  thus  there is $\delta_1=\delta_1(\varepsilon)>0$ such that  $0\leq x^1(t,0,\phi^1)\leq \varepsilon v^1$ for $t\ge 0$ whenever $0\leq \phi^1\leq \delta_1 v^1$.
Moreover, we find  $\widetilde\delta_1=\widetilde\delta_1(\widetilde\varepsilon)=\widetilde\delta_1(\varepsilon)>0$ such that if $0\leq \phi^1\leq \widetilde\delta_1 v^1$, then   $0\leq A_{21} x^1(t,0,\phi^1)\leq \widetilde\varepsilon v^2$ for $t\ge 0$.
\par
Take $\delta=\min(\delta_1,\widetilde\delta_1)$, and consider an initial condition $\phi\in C^+$ with $0\leq \phi^1\leq \delta v^1$ and $0\le \phi^2\le \varepsilon v^2$. We first solve the decoupled $n_1$-dimensional  system, and replace  in \eqref{4.8} the terms  $\sum_{p\in I_1}  a_{ip}\,x_p(t)$ by $\sum_{p\in I_1}  a_{ip}\,x_p(t,0,\phi^1)$. The crucial point is to check that $\varepsilon v^2$ is an ``upper" solution for this new  cooperative system, or in other words, that  $\widetilde f(t,\varepsilon v^2)\leq 0$ for any $t\geq 0$; this allows  concluding  that the set $[0,\varepsilon v^2]\subset C([-\tau,0];\R^{n_2})$ is positively invariant for \eqref{4.8} (see Lemma \ref{lem3}).
For each $i\in I_2$, we have
\begin{equation*}
\begin{split}
\widetilde f_i(t,\varepsilon v^2)&=-d_i\,\varepsilon v_i^2 +\dps \sum_{j\in I_2,j\not= i}  a_{ij}\,\varepsilon v_j^2 +  \beta_{i}\,\widetilde h_i(\varepsilon v_i^2)+\sum_{j\in I_1}  a_{ij}\,x^1_j(t,0,\phi^1)\\
&=-d_i\,\varepsilon v_i^2 +\dps \sum_{j\in I_2,j\not= i}  a_{ij}\,\varepsilon v_j^2 +\beta_i\varepsilon v_i^2  +  \beta_{i}(\widetilde h_i(\varepsilon v_i^2)-\varepsilon v_i^2)+\sum_{j\in I_1}  a_{ij}\,x^1_j(t,0,\phi^1)\,,
\end{split}
\end{equation*}
hence,
$$\widetilde f(t,\varepsilon v^2)\leq  \varepsilon M_{22}\, v^2 -\widetilde\varepsilon\, v^2+\widetilde\varepsilon\, v^2\leq 0.$$
As a consequence, and summarizing, we deduce that,  whenever $0\leq \phi^1\leq \delta v^1$ and $0\leq \phi^2\leq \varepsilon v^2$   then $0\leq x^1(t,0,\phi^1)\leq \varepsilon v^1$ and  $0\leq x^2(t,0,\phi^1,\phi^2)\leq y(t,0,\phi^1,\phi^2)\leq \varepsilon v^2$ for $t\ge 0$. This ends the proof.
\end{proof}

Theorem  \ref{thmOS} also provides conditions for partial extinction and partial persistence. As an illustration, we summarise the results for a Nicholson system.

\begin{exmp} {\rm Consider the Nicholson system with autonomous coefficients and time-dependent delays given by
\begin{equation}\label{3.16}
x_i'(t)=-d_ix_i(t)+\sum_{j=1,j\ne i}^n a_{ij}x_j(t)+\sum_{k=1}^m \be_{ik}    x_i(t-\tau_{ik}(t))e^{-c_{ik}x_i(t-\tau_{ik}(t))},\
i=1,\dots,n,\ t\ge 0,\end{equation}
where $d_i>0, c_{ik}>0, a_{ij}\ge 0, \be_{ik}\ge 0$ with $\be_i:=\sum_{k=1}^m \be_{ik} >0$, $\tau_{ik}:[0,\infty)\to [0,\tau]\ (\tau>0)$  are continuous, for all $i,j,k$, and $D-A$ is a non-singular M-matrix.  By applying Theorem  \ref{thmOS}  to this model, we obtain:
\begin{itemize}
\item[(i)] if $s(M)\le 0$, 0 is GAS;
\item[(ii)]  if $M$ is written in the triangular form (for some $k\in \{1,\dots, n\}$ and some $l\in \{1,\dots, k\}$)
 \begin{equation}\label{4.9}
M = \left[\begin{array}{cccccc}
M_{11} & \dots &0&0& \ldots & 0 \\
\vdots & \ddots  &\vdots & \vdots & & \vdots\\
0 & \dots &M_{ll}&0& \ldots& 0 \\
M_{l+1,1} & \dots & M_{l+1,l} &M_{l+1,l+1}& \ldots& 0 \\
\vdots &   &\vdots & \vdots &\ddots & \vdots\\
M_{k1} & \dots &  M_{k,l}& M_{k,l+1}  & \ldots& M_{kk}
\end{array}\right]\,,
\end{equation}
with $M_{jj} \ (1\le j\le k)$  irreducible blocks  and $M_{jp}\ne 0$ for some $p<j$ and $j=l+1,\dots k$, then \eqref{3.16}  is   permanent
 if and only if $s(M_{jj})>0$ for $j=1,\dots, l$;
 \item[(iii)] moreover, for $M$ written in the triangular form \eqref{4.9}, if $l>1$ and there exist $p,j\in \{ 1,\dots,l\}$  such that $s(M_{pp})\le 0$ and $s(M_{jj})>0$, then the $n_p$ populations $x_i(t)$ with $i\in I_p$ become extinct, whereas the $n_j$ populations $x_i(t)$ with $i\in I_j$ uniformly persist.
 \end{itemize}}
\end{exmp}

\begin{rmk}\label{rmk4.1} {\rm In this way,  we have recovered and extended all the results regarding extinction and uniform persistence established in  \cite{Faria11, FariaRost}
for the particular case of  \eqref{3.16} with constant delays $\tau_{ik}$ and $c_{ik}=1$ for all $i,k$.  For such autonomous systems, the sharp criterion for extinction of all populations, $s(M)\le 0$,  was proven in \cite{FariaRost} by using  the unimodal shape of the  specific Ricker nonlinearity $h(x)=xe^{-x}$.  However, as shown in the proof of Theorem \ref{thmOS},  the techniques presented  in Section 3, based on comparison of  solutions with  solutions for auxiliary cooperative systems,  allow us to    carry out the arguments for the more general model \eqref{3.14}. Also,  the permanence of Nicholson autonomous systems in \cite{FariaRost} was proven under the stronger requirement of $Mv\gg 0$ for some $v\gg 0$ (and $D-A$ a non-singular M-matrix).}
 \end{rmk}

\section*{Acknowledgements}
This work was partially supported by Funda\c c\~ao para a Ci\^encia e a Tecnologia under project UID/MAT/\-04561/2013 (T. Faria) and  by
Ministerio de Econom\'\i a y Competitividad under project MTM2015\--66330, and
the European Commission under project H2020-MSCA-ITN-2014 (R. Obaya and A.~M. Sanz).\\
The authors are very grateful to the referee, whose careful reading and valuable comments led to significant improvements of the manuscript.

\end{document}